\documentclass[10pt]{article}
\usepackage{amsmath, amsfonts, amsthm}
\usepackage{amssymb,mathrsfs, verbatim}
\usepackage{graphicx,amstext}
\usepackage{amsmath}
\usepackage{mathtools}
\usepackage{amssymb}
\usepackage{latexsym}
\usepackage{ mathrsfs }
\usepackage{enumerate}

\usepackage[latin1]{inputenc}
\newcommand{\clg}[1]{{\mathcal{#1}}}

\newcommand{\R}{\mathbb R}

\newcommand{\ve}{\varepsilon}
\newcommand{\vp}{\varphi}

\newcommand{\sgn}{\text{sgn}}

\newcommand{\ess}{\,{\rm ess}}

\newcommand{\dive}{{\rm div}}

\newcommand{\dist}{\,{\rm dist}}

\numberwithin{equation}{section}

\newtheorem{theorem}{Theorem}[section]
\newtheorem{proposition}{Proposition}[section]
\newtheorem{remark}{Remark}[section]
\newtheorem{lemma}{Lemma}[section]

\newtheorem{definition}{Definition}[section]

\begin{document}
\title{Dirichlet Problem for Degenerate Fractional 
\\
Parabolic Hyperbolic Equations}

\author{Gerardo Huaroto$^1$, Wladimir Neves$^2$}

\date{}

\maketitle

\footnotetext[1]{ Instituto de Matem\'atica, Universidade Federal
de Alagoas, Cidade Universit\'aria, Alagoas, Brazil  E-mail: {\sl gerardo.cardenas@im.ufal.br}}.\footnotetext[2]{
 Instituto de Matem\'atica, Universidade Federal
do Rio de Janeiro, C.P. 68530, Cidade Universit\'aria 21945-970,
Rio de Janeiro, Brazil. E-mail: {\sl wladimir@im.ufrj.br}.

}

\begin{abstract}
We are concerned in this paper with the 
degenerate fractional diffusion advection
equations posed in bounded domains. 
Due to a suitable formulation, 
we show the existence of weak entropy solutions
for measurable and bounded initial and Dirichlet boundary data. 
Moreover, we prove a $L^1-$type contraction property for 
weak entropy solutions obtained via parabolic 
perturbation. This is a weak selection principle which means that
the weak entropy solutions are stable in this class. 
\end{abstract}

\maketitle

\tableofcontents

\section{Introduction}
\label{INTRO}

We are concerned with the theory of evolutionary 
fractional partial differential equations
posed in bounded domains. 
Contrary to the $\R^n$ case,
there exist different non-equivalent definitions of the 
fractional Laplacian operator in bounded domains.  
We have this ambiguity not only for each 
particular type of the boundary condition considered, that is to say, Dirichlet, Neumann, etc., but also 
related to the process at the boundary, e.g. reflection, absorption, 
conditioning. Moreover, 
each particular problem may lead to a specific choice of the fractional Laplacian. 
This is an important issue, which we would like to 
stress the reader from the beginning.

\medskip
We consider isotropic degenerate fractional diffusion advection
equations, with general Dirichlet boundary 
condition. Hence the concept of entropy solutions are
needed, and due to this we assume the Regional (also called Censored) Fractional Laplacian operator, 
denoted by $(-\Delta)^{s}_{\Omega}$, ($0<s<1$), for the diffusion process. 
The main reason to consider the Regional Fractional Laplacian 
is not just related to integration by parts formula, see Guan \cite{Guan1}. In fact, 
we need (during the entropy formulation) the important property that, 
the fractional Laplacian of constant functions be zero.
This property clearly follows for the usual fractional Laplacian defined in $\R^n$, 
but not necessarily in bounded domains.  We return to this discussion at the end of this section. 

\medskip
More precisely,  we study in this paper the well-posedness of solutions
for the following initial-boundary value problem  
\begin{equation}
\label{FTPME} 
     \left \{
       \begin{aligned}
       \partial_t u+\dive f(u) + (-\Delta)^s_{\Omega}A(u)&= 0 \quad \,\,\, \text{in $Q_T$},
       \\
       u|_{t= 0} &= u_0 \quad \text{in $\Omega$},
       \\
       u &= u_b \quad \text{on $\Gamma_T$},
\end{aligned}
\right.
\end{equation}
where $Q_T:= (0,T) \times \Omega$, for any real number $T> 0$, and
$\Omega \subset \R^n$  is a bounded open set having smooth ($C^2$)
boundary $\Gamma$.
Moreover, $f:\R\rightarrow\R^n$ is called the flux function and $A:\R\rightarrow\R$ 
 is a nondecreasing function.
The initial data $u_0 \in L^\infty(\Omega)$, and 
the boundary condition $u_b \in L^\infty(\Gamma_T)$, 
where $\Gamma_T= (0,T) \times \Gamma$. 
The specific regularity for the functions $f$ and $A$ 
will be state a posteriori, also some extension condition on $u_b$, 
see Hypothesis (H1)-(H3) in Section \ref{SOLVAB}. 
 
 \medskip
As it is well known, 
the type of monotonicity assumption on the function 
$A(\cdot)$ allows some degeneration zones for the state variable $u$.
The degeneration makes 
the problem much different from the non-degenerate fractional parabolic one ($A'$ 
is strictly separated from 0) and from the hyperbolic one ($A'=0$), since the
 two behaviors are mixed together in such a way that depends on the solution itself. 
Therefore, we have to keep in mind both features of the fractional parabolic and of the hyperbolic equations.
Moreover, the problem is posed in bounded domains, albeit more realistic from the physical 
point of view it is much more difficult to treat. 

\medskip
In a similar context, Cifani, Jakobsen in 
\cite{cifane} considered the Cauchy problem for 
an analogous equation to $(\ref{FTPME})_1$, that is, 
$$
    \partial_t u + \dive f(u) + (-\Delta)^s A(u)= 0,\quad \text{ in $(0,T) \times \R^n$}, 
$$
hence in this case the fractional Laplacian operator $(-\Delta)^s$ 
can be defined using Fourier Transform by
$$
    \widehat{(-\Delta)^sf}(\xi)=|\xi|^{2s}\hat{f}(\xi), \quad \text{$s \in (0,1)$}. 
$$
Under some conditions for $f$ and $A$ they proved well-posedness of 
entropy solutions for $u_0 \in L^\infty(\R^n) \cap L^1(\R^n) \cap BV(\R^n)$.  
In fact, that paper were the first one in this direction, that is, they studied
the degenerate fractional diffusion advection equations,
and concerning bounded domains with non-homogeneous Dirichlet boundary data, this problem has 
never studied before as far as the authors know. 

\medskip
 It is important to observe right now that, the solvability of \eqref{FTPME}
 is established via the parabolic perturbation method, and the uniqueness is 
 obtained in this class. Indeed, the doubling variable techniques by Kru\v zkov
seems to be not adaptable since the non-local operator
presents some asymmetries during the doubling process. 
Moreover, we are dealing with bounded domains, hence a change of 
variables is needed to perform the technique,  (see for 
instance \cite{WN1}, Section 2), which complicates even more. 
 In fact, we only show a $L^1-$ type contraction, 
 see Theorem \ref{MainThm}, which implies uniqueness 
 and the well-posedness is just established when
 the boundary data is sufficiently regular, that is $C^1(\Gamma)$ with respect to the 
 spatial variable, see Remark \ref{UsualContraction}. Therefore, we obtain in the class
 of weak entropy solutions obtained via parabolic perturbation the inequality \eqref{L1-contration}, 
 which must coincide in the general uniqueness case. Similarly, we have the same observation
 to establish the uniqueness result via the kinetic method. Although, we claim that if it is possible 
 to obtain uniqueness, it should be the right way. We leave this important question to future studies.

\medskip
The theory of evolutionary fractional partial differential equations has received considerable 
interest in the last years, in particular since the seminal paper from Caffarelli, Vazquez \cite{Caffa}. 
Indeed, there exist a considerable amount of papers concerning anomalous diffusion, 
to mention a few \cite{NBUKGV, MSireVazquez, MBJLV, GHWN1,GHWN2, GKMKMK}.  
We should mention that, the equation \eqref{FTPME}
encompass various physical models, and so many applications considering anomalous diffusion in 
bounded domais (general Dirichlet data) could be represented by it.  
Moreover, this equation is a natural generalization to fractional diffusion setting of the  
degenerate diffusion advection equations  
\begin{equation}
\label{EQ100}
   \partial_t u+\dive f(u)+(-\Delta) A(u)= 0 \quad \text{posed in bounded domains.}
\end{equation}

One remarks that, there exists a list of important papers related to \eqref{EQ100}. 
Indeed, Carrillo \cite{Carrillo} has proved existence and uniqueness of entropy solutions for homogeneous 
boundary condition, $(u_b= 0)$. Mascia, Porreta, Terracina \cite{Mascia}, and Michel, Vovelle \cite{Vovelle} extended that
result for non-homogeneous boundary conditions. Karlsen, Risebro \cite{KarlsenRisebro} considered the case of rough coefficients. 
Moreover, we would like to mention the anisotropic diffusion
correlated results. For instance, we address the reader to Bendahmane, Karlsen  \cite{BendahmaneKarlsen}, 
and Kobayasi, Ohwa \cite{KobayasiOhwa}. Clearly, the above list is not exhaustive.

\medskip
Let us return to the assumption of the Regional Fractional Laplacian. As mentioned at the begging,
it is well known in the literature that the fractional Laplacian $(-\Delta)^s$ has many equivalent, albeit formally different
definitions. For instance, the fractional Laplacian
can also be described using singular integrals in the following way
\begin{equation}
\label{DEFFRACRN}
(-\Delta)^s u(x)= C_{n,s} \, \lim_{\epsilon \to 0^+} \! \int_{\mathbb
R^n \setminus B_\epsilon(x)}\frac{u(x)-u(y)}{|x-y|^{n+2s}} \ dy,
\end{equation}
where $B_\epsilon(x)$ is the closed ball centered 
at $x \in \R^n$ with radius $\epsilon> 0$, and 
$$
   C_{n,s}=  \frac{\Gamma(\frac{n}{2} +s)}{\pi^{2 s + \frac{n}{2}} \Gamma(-s)}. 
$$
It follows from \eqref{DEFFRACRN} that $(-\Delta)^s u(x)= 0$ when $u$ is a constant function,
which is a required feature. Albeit, it is well known (see \cite{HohJacob}) that the Dirichlet 
type problem  
\begin{equation}
\label{DRFL}
   (-\Delta)^s u= 0 \quad \text{in $\Omega$}, \quad u= u_b \quad \text{on $\Gamma$}, 
\end{equation}
is not well posed. In order to have well-posedness we need to assume $u= u_b$ on $\R^n \setminus \Omega$,
which is not suitable since we would like to set locally the boundary data at $\Gamma$.  
To follow, let us recall from \cite{GHWN1} the Dirichlet spectral fractional Laplacian operator, denoted by $(-\Delta_D)^s$, 
which can be defined as a fractional power 
of the Laplacian operator $(-\Delta_D)$ with homogeneous Dirichlet boundary condition, or equivalently by 
$$
\begin{aligned}
(-\Delta_D)^s u(x)&= \frac{1}{\Gamma(-s)} \int_0^\infty \big( \int_\Omega K_\Omega(t,x,y) u(y) dy - u(x) \big) \frac{dt}{t^{1+s}}
\\[5pt] 
&= \frac{1}{\Gamma(-s)} \int_0^\infty \!\!  \int_\Omega K_\Omega(t,x,y) \big(u(y) - u(x) \big) dy \frac{dt}{t^{1+s}}
\\[5pt] 
& + \frac{u(x)}{\Gamma(-s)} \int_0^\infty \!\! \int_\Omega \big(K_\Omega(t,x,y) - 1 \big) dy \frac{dt}{t^{1+s}}, 
\end{aligned}
$$
where $K_\Omega$ is the Heat Kernel associated to $(-\Delta_D)$. 
Therefore, it does not happen that, $(-\Delta_D)^s u(x)= 0$ when $u(x)$ is a constant function.
Hence even if we are allowed to specify the boundary condition locally for $(-\Delta_D)^s$,
this is not the correct choice.   

\medskip
On the other hand, if $u$ is a constant function, then from the definition of the 
Regional Fractional Laplacian, we have $(-\Delta)^{s}_{\Omega} u(x)= 0$. Moreover, it has been shown by 
Guan, Ma \cite{GuanMa} that, the Dirichlet type problem for the Regional Fractional Laplacian \eqref{DRFL}
is well-posed, in particular for $s \in (1/2, 1)$ when the censored process could approach to the boundary.  
Consequently, we assume in this paper the  $(-\Delta)^s_\Omega$ nonlocal operator to perform the diffusion process. 

\subsection{Functional spaces}

Let $U$ be an open set in $\R^n$. We denote by $L^p(U)$ the set of real $p-$summable functions 
with respect to the Lebesgue measure (vector ones should be understood 
componentwise). 

$\bullet$ {\bf The space $W^{s,p}(U)$}

The Sobolev space is denoted by $W^{s,p}(U)$, where a 
real $s\geqslant 0$ is the
smoothness index, and a real $p\geqslant 1$
is the integrability index. More precisely, for $s \in (0,1)$, 
$p \in [1,+\infty)$, the fractional Sobolev space of
order $s$ with Lebesgue exponent $p$ is defined by
$$
   W^{s,p}(U):= \Big\{ u \in L^p(U): \int_{U} \int_{U}\dfrac{\vert u(x)-u(y)\vert^p}{
\vert x-y\vert ^{n+sp}} \ dx dy< + \infty \Big\},
$$
endowed with norm
$$
   \Vert u\Vert_{W^{s,p}(U)}= \left(  \int_{U}\vert u\vert^p dx + \int_{U} \int_{U}\dfrac{\vert u(x)-u(y)\vert^p}{
    \vert x-y\vert ^{n+sp}}dxdy \right)^\frac{1}{p}.
$$
Moreover, for $s > 1$ we write $s = m + \vartheta$, where $m$ is an integer
and $\vartheta \in (0, 1)$. In this case, the space $W^{s,p}(U)$ consists of those equivalence classes
of functions $u \in W^{m,p}(U)$ whose distributional derivatives $D^{\alpha} u$, with $|\alpha| = m$,
belong to $W^{\vartheta,p}(U)$, that is
$$
        W^{s,p}(U)= \Big\{ u \in W^{m,p}(U): 
        {\displaystyle \sum_{\vert \alpha\vert = m}}\Vert D^{\alpha}u \Vert_{W^{\vartheta,p}(U)}< \infty \Big\},
$$
which is a Banach space with respect to the norm
$$
\Vert u\Vert_{W^{s,p}(U)}= \Big(  \Vert u\Vert^p_{W^{m,p}(U)} 
+ {\displaystyle \sum_{\vert \alpha\vert = m}}\Vert D^{\alpha}u \Vert^p_{W^{\vartheta,p}(U)} \Big)^\frac{1}{p}.
$$
If $s = m$ is an integer, then the space $W^{s,p}(U)$ coincides with the Sobolev space
$W^{m,p}(U)$. Also, it is very interesting the case when $p = 2$, i.e. $W^{s,2}(U)$. 
In this case, the (fractional)
Sobolev space is also a Hilbert space, and we can consider the inner product
\begin{equation}
\label{G_inner}
    \langle u,v\rangle_{W^{s,2}(U)}= \langle u,v\rangle + \int_{U} \int_{U}
    \frac{(u(x)-u(y))}{\vert x-y\vert ^{\frac{n}{2}+s}} \ \frac{(v(x)-v(y))}{\vert x-y\vert ^{\frac{n}{2}+s}} \ dx dy,
\end{equation}
where $\langle\cdot,\cdot\rangle$ is the inner product in $L^2(U)$, and the term 
\begin{equation}
\label{G_semi_inner_prod}
    [u,v]_{W^{s,2}(U)}= \iint_{U \times U} 
    \frac{(u(x)-u(y))(v(x)-v(y))}{\vert x-y\vert ^{n+2s}} \ \ dx dy,
\end{equation}
is  Gagliardo semi-inner product.
Moreover, we can define the subspace $W_0^{s,p}(U)$ as
$$
W^{s,p}_0(U)= \overline{C^{\infty}_c(U)}^{\|\cdot\|_{W^{s,p}(U)}}.
$$

$\bullet$ {\bf The space $H^s(U)$, 
$H^s_0(U)$ and $H^{1/2}_{00}(U)$}

We can define for $s\in (0,1)$, (see Lions, Magenes \cite{LionsMagenes}), the spaces $H^s(U)$ 
by interpolation between $H^1(U)$ and $L^2(U)$, that is 
$$
H^s(U)= [H^1(U),L^2(U)]_{1-s}.
$$
According to this definition, this space is a Hilbert space with the natural norm given by the interpolation.
One remarks that, if $U$ has Lipschitz boundary the spaces $W^{s,2}(U)$ and $H^s(U)$ are
equivalent, see Tartar $\cite{Tartar}$ p. 83, and p. 84.

As mentioned above, the space $H^s_0(U)$ can be defined by
$$
H^s_0(U)= \overline{C^{\infty}_c(U)}^{\|\cdot\|_{H^s(U)}}, 
$$
and once $U$ has Lipschitz boundary,
there exists an equivalent definition given via interpolation, namely Theorem 11.6 of \cite{LionsMagenes}, which
states that
$$
H^s_0(U)=[H^1_0(U),L^2(U)]_{1-s},
$$
for each $s\in (0,1) \setminus \{1/2\}$. The case $s = 1/2$ is special and generates the so called 
Lions-Magenes space $H^{1/2}_{00}(U)$, which is defined by
$$
    H^{1/2}_{00}(U):= [H^1_0(U),L^2(U)]_{1/2},
$$
and has the following characterization
$$
H^{1/2}_{00}(U)= \Big\{ u\in H^{1/2}(U);\int_{U}\dfrac{u(x)^2}{\dist(x,\partial U)}dx<\infty \Big\}. 
$$

Finally, we recall the following results.
\begin{theorem}
Let $\Omega$ be a bounded open set in $\R^n$ with Lipschitz boundary. Then, $C^{\infty}_c(\Omega)$ is
dense in $H^s(\Omega)$ if, and only if, $0 < s \leq 1/2$, and in this case,  
$H^s_0(\Omega) = H^s(\Omega)$. If $s > 1/2$, then $H^s_0(\Omega)\subset H^s(\Omega)$ and the inclusion is strict.
\end{theorem}
\begin{proof}
See Lions, Magenes \cite{LionsMagenes} Theorem 11.1, and Tartar \cite{Tartar}, p. 160.
\end{proof}

\begin{theorem}
\label{TraceTHM}
Let $\Omega$ be a bounded open set in $\R^n$ with Lipschitz boundary. 
For $1/2 < s \leq 1$, there exists a surjective linear operator 
$T: H^s(\Omega) \rightarrow H^{s - \frac{1}{2}}(\partial\Omega)$, such that,
$T^{-1}(0)= H^s_0(\Omega)$, and 
$$
   H_0^s(\Omega) \equiv \{ u \in H^s(\Omega): \text{$u= 0$ on $\partial\Omega$ in the sense of trace}\}. 
$$
\end{theorem}
\begin{proof}
See Lions, Magenes \cite{LionsMagenes} Theorem 9.4, and Theorem 11.5.
\end{proof}

\section{Regional Fractional Laplacian} 
\label{REGFRACLAP}

In this section we consider the background of the Regional Fractional Laplacian operator
used in this paper. In fact, we mainly provide the proofs of the new results, and 
address the reader to \cite{Guan1, Guan} 
for an introduction of this subject.  

\medskip
Let $\Omega \subset \R^n$ be an open set, $s \in (0,1)$, and for $x \in \Omega$ we define the following Lebesgue-Stieltjes measures on $\Omega$ 
$$
   d \mu_s(x):= \frac{1}{(1+|x|)^{n+2s}} \ dx. 
$$
Hence a measurable function $f  \in L^1(\Omega; d\mu_s)$,  means that 
$$
    \|f\|_{L^1_\mu}:=  \int_\Omega |f(x)| \ d\mu_s(x)< \infty.
$$
Now, we define for $u \in L^1(\Omega; d\mu_s)$, any $\epsilon> 0$ fixed, and almost every $x \in \Omega$ 
$$
   (-\Delta)^{s}_{\Omega,\epsilon}u(x):= C_{n,s}
   \int_{\Omega \setminus B_\epsilon(x)} \frac{u(x)-u(y)}{|x-y|^{n+2s}} \ dy.
$$
Clearly, for each $\epsilon> 0$ the operator $(-\Delta)^{s}_{\Omega,\epsilon}$ is well-defined, 
non-local and depends on $\Omega$. 
Moreover, if the function $u(x)$ is constant, then $(-\Delta)^{s}_{\Omega,\epsilon} u(x)= 0$ for each $\epsilon> 0$. 

\begin{definition}
Given $u \in L^1(\Omega; d\mu_s)$ and $0<s <1$, the Regional Fractional Laplacian, denoted 
$(-\Delta)_{\Omega}^{s}$, is defined as
$$
(-\Delta)^{s}_{\Omega}u(x):= \lim_{\epsilon\to0^+}(-\Delta)^{s}_{\Omega,\epsilon}u(x)
$$
for almost all $x \in \Omega$, provided the above limit exists.
\end{definition}

One remarks that, if $\Omega= \R^n$ then the 
regional fractional Laplacian $(-\Delta)_{\R^n}^{s}$ coincides with the usual $s-$fractional Laplacian on $\R^n$. 
Furthermore, all the properties mentioned above to $(-\Delta)^{s}_{\Omega,\epsilon}$ remains 
valid when we pass to the limit as $\epsilon \to 0$, when it exists. 

\medskip
The following result is a useful product formula, obtained in \cite{Guan} Lemma 3.5, for the Regional Fractional Laplacian.
\begin{lemma}
\label{LFprodu}
Let $\Omega \subset \R^n$ be an open set, $u, v$ be two real value measurable functions such that $u v \in L^1(\Omega; d\mu_s)$. 
If $(-\Delta)^{s}_\Omega u(x)$, $(-\Delta)^{s}_\Omega v(x)$ exist for almost all $x \in \Omega$, and also
$$
   \int_\Omega \frac{|(u(x)-u(y))(v(x)-v(y))|}{|x-y|^{n+2s}}dy< \infty,
$$   
then for almost all $x \in \Omega$,
$$
\begin{aligned}
(-\Delta)^{s}_\Omega(u v)(x)=& v(x) \ (-\Delta)^{s}_\Omega u(x) + u(x) \ (-\Delta)^{s}_\Omega v(x)
\\[5pt]
&- C_{n,s} \int_\Omega\frac{(u(x)-u(y))(v(x)-v(y))}{|x-y|^{n+2s}} dy.
\end{aligned}
$$
\end{lemma}

Now we present some definitions and notations,
which are necessary to state the integration by parts formula 
for the Regional Fractional Laplacian. 

\medskip
Henceforth $\Omega$ is a bounded open set in $\R^n$ with a $C^2-$ boundary $\Gamma$, and assume $\frac{1}{2}< s \leq 1$. 
Then, we can find $\delta> 0$ depending on $\Omega$, and a
function $h \in C^2(\Omega)$ depending also on $\Omega$, such that, for $\delta> 0$ sufficiently small
$$
    h(x)= d(x)^{2s-1} \quad \text{for any $x \in \Omega_\delta$}, 
$$
where $d(x)= \dist(x,\partial\Omega) \big)$ and 
\begin{equation}
\label{Omega_delta}
  \Omega_\delta=\{x \in \Omega ; 0< d(x)< \delta\}. 
\end{equation}
Conveniently, we define for some $f, g \in C^2(\overline{\Omega})$
$$
C^2_{2s}(\overline{\Omega}):= \left\lbrace u:\,u(x)=f(x)h(x)+g(x), \forall x \in \Omega \right\rbrace.
$$

\begin{remark}
The set $C^2_{2s}(\overline{\Omega})$ does not depend on the choice of $h$. Moreover, we may always assume that,
functions in $C^2_{2s}(\overline{\Omega})$ is 
defined on $\overline{\Omega}$ by continuous extension. 
Precisely, we have 
$$
    C^2(\overline{\Omega})\subset C^2_{2s}(\overline{\Omega})\subset C(\overline{\Omega})\cap C^2(\Omega).
$$ 
In particular, for $s= 1$ it follows that $C^2_{2}(\overline{\Omega})= C^2(\overline{\Omega})$.
\end{remark}

The next lemma is important to show that, $u \in H^s(\Omega)$ whenever 
 $u\in C^2_{2s}(\overline{\Omega})$. More precisely, we have the following
\begin{lemma}
\label{Lemminclu}
Let $\frac{1}{2}< s< 1$. If $u\in C^2_{2s}(\overline{\Omega})$, then $u\in H^{s}(\Omega)$.
\end{lemma}
\begin{proof}
See Lemma 3.1 in \cite{Guan1}.
\end{proof}
Analogously, we have the following
\begin{proposition}
Given $u \in C^2_{2s}(\overline{\Omega})$ for any $\frac{1}{2} < s < 1$, then the Regional Fractional Laplacian
of the function $u$ exists, and $(-\Delta)^{s}_\Omega u \in L^2(\Omega)$.
\end{proposition}
\begin{proof}
1. First, the existence of $(-\Delta)^{s}_\Omega u$ follows directly from Proposition 2.2 in \cite{Guan}.
Recall that $u \in C^2_{2s}(\overline{\Omega}) \subset C(\overline{\Omega})\cap C^2(\Omega)$. 
Similarly, it follows from item $(ii)$ in Proposition 2.3 in \cite{Guan} that, 
$(-\Delta)^{s}_\Omega u \in L^2(\Omega)$. 
\end{proof}

Following Guan \cite{Guan1}, with minor modifications, 
the following definition generalizes the concept of normal derivative of the function $u$ in direction 
of the outer unitary normal vector field $\nu$ to $\Omega$, that is, 
$\nabla u \cdot \nu \equiv \partial_\nu u$ on $\Gamma$.
\begin{definition}
For $\sigma= 2s -1$, $\frac{1}{2} < s < 1$, 
the fractional normal derivative 
of the function $u$, denoted $\partial^{\sigma}_\nu u$, 
is defined for $r \in \Gamma$ and $\tau> 0$ by 
$$
   \partial^{\sigma}_\nu u(r):=  \lim_{ \tau \to 0^+} \big(\sigma \dfrac{u(r) - u(r - \tau \nu(r))}{\tau^{\sigma}} \big),
$$
provided this limit exists. 
\end{definition}
Then, we have the following 
\begin{lemma}
Given $u \in C^2_{2s}(\overline{\Omega})$ for any $\frac{1}{2} < s < 1$, the fractional normal derivative 
of the function $u$ exists, and $\partial^{\sigma}_\nu u(r)= -\sigma \, f(r)$ for each $r \in \Gamma$. 
\end{lemma}

\begin{proof}
Since $u \in C^2_{2s}(\overline{\Omega})$, for $\tau> 0$ sufficiently small, it follows by definition that 
$$
\begin{aligned}
   \sigma \dfrac{u(r) - u(r - \tau \nu(r))}{\tau^{\sigma}}&= \sigma \dfrac{f(r) h(r) + g(r)}{\tau^{\sigma}} - \sigma
   \dfrac{f(r - \tau \nu(r)) \tau^\sigma + g(r - \tau \nu(r))}{\tau^{\sigma}}
   \\[5pt]
   &= - \sigma \, f(r - \tau \nu(r)) + \sigma \dfrac{g(r) - g(r - \tau \nu(r))}{\tau} \, \tau^{2-2s},
\end{aligned}   
$$
where we have used that $h(r)= 0$ and $\sigma= 2s - 1$. 
Then, passing to the limit as $\tau \to 0^+$, the lemma is proved. 
\end{proof}

\begin{remark}
\label{FracTrace}
One remarks that, if $u\in C^1(\overline{\Omega})$, then $\partial^\sigma_\nu u(r)= 0,$ for all $r \in \Gamma$,
that is to say, the fractional normal derivative 
of a sufficiently regular (up to the boundary) function $u$ is zero. 
Moreover, from the proof of the above lemma, it is enough to have $f, g \in C(\overline{\Omega})$, with 
$$
    \dfrac{g(r) - g(r - \tau \nu(r))}{\tau} \quad \text{uniformly bounded in $\overline{\Omega}_\delta$}
$$
to ensure the existence of the fractional normal derivative 
of a function $u$, such that, $u(x)= f(x) \, h(x) + g(x)$. 
\end{remark}

At this point, we are able to state plainly the 
following integration 
by parts formula for the Fractional Regional Laplacian, see Theorem 3.3 in \cite{Guan1}. 

\begin{theorem}
\label{Thintpart}
Let $\frac{1}{2}<s<1$, $u\in C^{2}_{2s}(\overline{\Omega})$, $v\in H^{s}(\Omega)$ and $\sigma= 2s-1$. 
Then, the following integration by parts formula holds 
\begin{equation}
\begin{aligned}
\int_{\Omega}v(x)(-\Delta)^s_\Omega u(x) \ dx
=&\frac{C_{n,s}}{2}\iint_{\Omega\times\Omega}\dfrac{(v(x)-v(y))(u(x)-u(y))}{|x-y|^{n+2s}}dx\,dy
\\
&-\mathcal{N}_\sigma \int_{\Gamma}v(r) \ \partial^{\sigma}_\nu u(r) \ dr, 
\end{aligned}\label{eq:intpart}
\end{equation}
where 
$$ 
   \mathcal{N}_\sigma= \frac{C_{1,(\sigma+1)/2}}{(\sigma+1) \sigma} 
   \int_0^\infty \Big( \frac{1}{|\tau - 1|^\sigma} - \frac{1}{(\max\{\tau,1\})^\sigma} \Big) \tau^{\sigma-1} \ d\tau. 
$$ 
\end{theorem}
Comparing \eqref{eq:intpart} and the classical Green Formula for the Laplace operator, 
we observe that $\mathcal{N}_\sigma \partial^\sigma_\nu$  plays the role (for the Regional 
Fractional Laplace operator) that $\partial_\nu u$ does for 
the Laplace operator. 



\medskip 
Thanks for the above results, we extend to functions in $W^{s,2}(\Omega)$
the definition of $(-\Delta)^s_\Omega$, for $s \in(\frac{1}{2},1)$.  To begin, 
we take $u \in C^{2}_{2s}(\overline{\Omega})$ and from \eqref{Thintpart}, 
it follows for each $v \in W^{s,2}_0(\Omega)$ that 
$$
 \int_{\Omega}v(x)(-\Delta)^s_\Omega u(x) \ dx
=\frac{C_{n,s}}{2}\iint_{\Omega\times\Omega}\dfrac{(v(x)-v(y))(u(x)-u(y))}{|x-y|^{n+2s}}dx\,dy.
$$
Therefore, we have 
$$
\begin{aligned}
\big| \int_{\Omega}v(x) \ (-\Delta)^s_\Omega u(x) dx \big|&
\leq \frac{C_{n,s}}{2}\iint_{\Omega\times\Omega}\dfrac{|v(x)-v(y)| \ |u(x)-u(y)|}{|x-y|^{n+2s}}dx\,dy
\\[5pt]
&\leq C  \|u\|_{W^{s,2}(\Omega)} \,\|v\|_{W^{s,2}(\Omega)},
\end{aligned}
$$
where $C$ is a positive constant and we have used \eqref{G_inner}, which implies that 
\begin{equation}
\|(-\Delta)^s_\Omega u\|_{W^{-s,2}(\Omega)}\leq C \|u\|_{W^{s,2}(\Omega)}.
\label{eq:ext1}
\end{equation}
Now, we assume $u\in W^{s,2}({\Omega})$ and by a simple density argument, 
there exits a sequence $\{u_n\}_{n=1}^{\infty}$, 
$u_n \in C^{2}_{2s}(\overline{\Omega})$ for $n \geq 1$, such that $u_n$ converges to
$u$ in $W^{s,2}(\Omega)$, as $n \to \infty$. According to \eqref{eq:ext1}, we may write 
\begin{equation}
\|(-\Delta)^s_\Omega u_n-(-\Delta)^s_\Omega u_m\|_{W^{-s,2}(\Omega)}\leq C \|u_n-u_m\|_{W^{s,2}(\Omega)}.
\label{eq:ext2}
\end{equation}
Therefore, $\{(-\Delta)^s_\Omega u_n\}_{n=1}^{\infty}$ is a Cauchy sequence in $W^{-s,2}(\Omega)$, 
which is a Hilbert space. 
Then, we may define for each $u \in W^{s,2}(\Omega)$,
$$
    (-\Delta)^s_{\Omega} u:=\lim_{n\to 0} (-\Delta)^s_\Omega u_n\quad\mbox{ in } W^{-s,2}(\Omega).
$$
Due to \eqref{eq:ext2} this definition does not depend on the choice of the approaching sequence
to the function $u$. Then, we have proved the following 
\begin{proposition}
\label{cor3part}
For each $\frac{1}{2}< s< 1$, $u \in W^{s,2}({\Omega})$, and $v \in W^{s,2}_0({\Omega})$, it follows that 
\begin{equation}
\label{EXT1}
     \left\langle (-\Delta)^s_{\Omega} u , v\right\rangle =\frac{C_{n,s}}{2}\iint_{\Omega\times\Omega}\dfrac{(v(x)-v(y))(u(x)-u(y))}{|x-y|^{n+2s}}dx\,dy,
\end{equation}
where $\langle\cdot,\cdot\rangle$ denotes the pairing between $W^{-s,2}(\Omega)$ and $W^{s,2}_0(\Omega)$.
\end{proposition}

\subsection{General Green Formula}

Motivated by the Proposition \ref{cor3part}, and since the right hand side of \eqref{EXT1} is well-defined for any $s \in (0,1)$, 
we extend the definition of the Regional Fractional Laplacian to functions in $W^{s,2}(\Omega)$, for all $s\in(0,1)$. Thus we consider the following 
\begin{definition}
Let $\Omega \subset \R^n$ be a bounded open set with $C^2$ boundary, and $s \in (0,1)$. Given a function $u \in W^{s,2}(\Omega)$,
the Regional Fractional Laplacian $(-\Delta)_\Omega^s u$ is defined as an element of $W^{-s,2}(\Omega)$, that is, 
for each $v\in W^{s,2}_0(\Omega)$
$$
\left\langle(-\Delta)^s_\Omega u , v\right\rangle:=
\frac{C_{n,s}}{2}\iint_{\Omega\times\Omega}\dfrac{(v(x)-v(y))(u(x)-u(y))}{|x-y|^{n+2s}}dx\,dy.
$$
Moreover,  we define
$$
D\big((-\Delta)_\Omega^{s}\big):=\left\lbrace u\in W^{s,2}(\Omega): (-\Delta)^{s}_{\Omega}u\in L^2(\Omega) \right\rbrace.
$$
\end{definition}

Clearly, for $u$ sufficiently regular the definition above coincides 
with the previous one, hence we maintain the notation i.e. $(-\Delta)_\Omega^s u$. 
Moreover, from the above definition if $u \in D\big((-\Delta)_\Omega^{s}\big)$, then
for each $v\in W^{s,2}_0(\Omega)$, 
$$
\int_{\Omega} v(x) \, (-\Delta)^s_\Omega u(x)\,dx=
\frac{C_{n,s}}{2}\iint_{\Omega\times\Omega}\dfrac{(v(x)-v(y))(u(x)-u(y))}{|x-y|^{n+2s}}dx\,dy.
$$
%

Now, we establish the trace operator $\mathcal{N}_\sigma \partial_\nu^\sigma u$, for functions 
$u \in D\big((-\Delta)_\Omega^s\big)$ and any $\sigma  \in (0,1)$, that is $\frac{1}{2}< s< 1$. 
To this end, we define for each function ${v} \in H^{s-\frac{1}{2}}(\Gamma)$, 
\begin{equation}
\label{tracefract}
\begin{aligned}
\langle \mathcal{N}_\sigma \partial^\sigma_\nu u,v \rangle_\Gamma :=
&\frac{C_{n,s}}{2}\iint_{\Omega\times\Omega}\dfrac{(\tilde{v}(x)-\tilde{v}(y))(u(x)-u(y))}{|x-y|^{n+2s}}dx\,dy
\\[5pt]
&- \int_{\Omega}  \tilde{v}(x)\,(-\Delta)^s_\Omega u(x) dx,
\end{aligned}
\end{equation}
where $\langle\cdot,\cdot\rangle_\Gamma$ denotes the inner product in $H^{s-\frac{1}{2}}(\Gamma)$, 
and $\tilde{v} \in W^{s,2}(\Omega)$ is an extension of $v$. 
Then, we have the following 
\begin{theorem} For any $s \in (\frac{1}{2}, 1)$, $\sigma= 2 s - 1$,   
the linear operator 
$$
 \mathcal{N}_\sigma \partial_\nu^\sigma: D\big( (-\Delta)_\Omega^s\big) \rightarrow H^{s-\frac{1}{2}}(\Gamma) $$
is well-defined and continuous.
\end{theorem}
\begin{proof}
1. First, we show that the definition of $\langle \mathcal{N}_\sigma \partial^\sigma_\nu u,v \rangle_\Gamma$ is 
independent of the extension $\tilde{v}$ of $v$. From Theorem \ref{TraceTHM}, 
let $\tilde{v}_1,\,\tilde{v}_2 \in W^{s,2}(\Omega)$ 
be such that, $T(\tilde{v}_1)= T(\tilde{v}_2)= v$. Therefore, we obtain from \eqref{tracefract}
$$
\begin{aligned}
&\langle \mathcal{N}_\sigma \partial^\sigma_\nu u, T(\tilde{v}_1-\tilde{v}_2)\rangle_\Gamma
\\[5pt]
&\qquad = \frac{C_{n,s}}{2}\iint_{\Omega\times\Omega}\!\!\!\dfrac{((\tilde{v}_1-\tilde{v}_2)(x)-(\tilde{v}_1-\tilde{v}_2)(y))(u(x)-u(y))}{|x-y|^{n+2s}}dx\,dy
\\[5pt]
&\qquad - \int_{\Omega} (\tilde{v}_1-\tilde{v}_2)(x) \, (-\Delta)^s_\Omega u(x)  \, dx. 
\end{aligned}
$$
Hence assuming that the right hand side of the above equality is different from zero,
we obtain a contradiction, since $\tilde{v}_1-\tilde{v}_2\in W^{s,2}_0(\Omega)$. 
 
\medskip
2. Now, we prove the boundedness of $\mathcal{N}_\sigma \partial_\nu^\sigma$, the linearity is trivial.  
Let $\mathcal{E}$ be the extension operator from $H^{s-\frac{1}{2}}(\Gamma)$ to $W^{s,2}(\Omega)$, thus
we have the following estimate
$$
\begin{aligned}
\| \mathcal{N}_\sigma \partial_\nu^\sigma u\|_{H^{s-\frac{1}{2}}(\Gamma)}&
= \!\!\!  \sup _{v \in H^{s-\frac{1}{2}}(\Gamma) \backslash\{0\}} \!\!\! \frac{|\langle\mathcal{N}_\sigma \partial_\nu^\sigma u, v\rangle_{\Gamma}|}{\|v\|_{L^2(\Gamma)}} 
\leq C \!\!\!  \sup _{v \in H^{s-\frac{1}{2}}(\Gamma) \backslash\{0\}} \!\! \frac{|\langle\mathcal{N}_\sigma \partial_\nu^\sigma u, v\rangle_{\Gamma}|}{\|\mathcal{E} v\|_{W^{s,2}(\Omega)}}
\\[5pt]
&\leq C \left(\,\|u\|_{W^{s,2}(\Omega)}+\|(-\Delta)^s_\Omega u\|_{L^2(\Omega)}\,\right),
\end{aligned}
$$
which implies that $\mathcal{N}_\sigma \partial_\nu^\sigma$ is a continuous operator. 
\end{proof}
The definition for $ \mathcal{N}_\sigma \partial^\sigma_\nu u$ over $\Gamma$ yields the 
General Green Formula for $(-\Delta)_\Omega^s$
\begin{equation}
\label{GGF}
\begin{aligned}
\int_{\Omega} \vp(x) \, (-\Delta)^s_\Omega u(x) \, dx&=
\frac{C_{n,s}}{2}\iint_{\Omega\times\Omega}\dfrac{(\vp(x)-\vp(y))(u(x)-u(y))}{|x-y|^{n+2s}}dx\,dy
\\[5pt]
& \quad - \mathcal{N}_\sigma \int_{\Gamma} \vp(r) \ \partial^{\sigma}_\nu u(r) \ dr, 
\end{aligned} 
\end{equation}
for any $\vp \in C^\infty_c(\R^n)$ and $u \in D\big( (-\Delta)_\Omega^s\big)$. 

\begin{remark}
\label{TRACECOND}
If $u \in H^{2s}(\Omega)$ for $s \in (1/2,1)$, then $(-\Delta)_\Omega^s u \in L^2(\Omega)$, 
(see for instance Sections 2.2 and 2.3 in $\cite{BGU}$), hence $u \in D\big( (-\Delta)_\Omega^s\big)$.

\end{remark}

%
%
\section{Entropy formulation} 
\label{SOLVAB}

The aim of this section is to introduce the entropy formulation for the isotropic 
degenerate fractional diffusion advection equations 
$$
  \partial_t u+\dive f(u) + (-\Delta)_{\Omega}^sA(u)= 0,
$$
posed in $Q_T= (0,T)\times\Omega$, for any $0<s<1$. 
Since we seek for suitable week solutions and 
the above equation degenerates, the function $u$ 
may admit shocks. Therefore, in order to establish an admissibility criteria 
to may select the correct physical solution, we need the concept of entropy. 

\begin{definition}
\label{DEFENT}
A pair $\mathbf{F}(u)= (\eta(u), q(u))$ is called an entropy pair for the first equation in \eqref{FTPME}, 
if there exists $\eta: \mathbb{R} \rightarrow \mathbb{R}$ a Lipschitz continuous and also convex 
function and the function $q: \mathbb{R} \rightarrow \mathbb{R}$, which satisfies for almost all $u \in \mathbb{R}$
$$
q^{\prime}(u)=\eta^{\prime}(u) f^{\prime}(u).
$$
We call $\eta(u)$ an entropy and $q(u)$ the associated entropy flux. 
Moreover, $\mathbf{F}(u)$ is called a generalized entropy pair when 
it is the uniform limit of a sequence of entropy pairs over compact sets. 
\end{definition}

The most important example of a generalized entropy pairs are 
the Kru\v zkov's entropies, namely 
$$
   \mathbf{F}(u, k)=(|u-k|, \operatorname{sgn}(u-k)(f(u)-f(k)))
$$
for each $k \in \R$.  Another two examples of generalized family of entropy pairs,
which will be conveniently used are the so 
called Kru\v zhkov semi-entropies, i.e. 
$$
\begin{aligned}
\mathbf{F}^{\pm}(u,k)&=\left(|u-k|^{\pm}, \operatorname{sgn}^{\pm}(u-k)(f(u)-f(k))\right)
\\
&=:(\eta^\pm_k(u), q^\pm_k(u))
\end{aligned}
$$
for each $k \in \mathbb{R}$, where $|v|^{\pm}:=\max \{\pm v, 0\}$, and
$$
\operatorname{sgn}^{+}(v):=\left\{\begin{array}{ll}
1, & \text { if } v>0 \\
0, & \text { if } v \leqslant 0,
\end{array} \quad \operatorname{sgn}^{-}(v):=\left\{\begin{array}{ll}
\!\! -1, & \text { if } v<0 \\
\; 0, & \text { if } v \geqslant 0.
\end{array}\right.\right.
$$

\medskip
At this point we make the following Hypothesis:
\begin{itemize}
\item[(H1)] The flux function $f \in C^1(\R)$ and $f'$ is locally Lipschitz continuous.

\item[(H2)] The nondecreasing function $A \in C^1(\R)$, $A^\prime$  is locally Lipschitz continuous, 
and without loss of generality $A(0)= 0$.

\item[(H3)] The boundary condition $u_b \in L^\infty(\Gamma_T)$ is the trace of a function 
$$
    \text{$ \tilde{u}_b \in L^{\infty}(Q_T)$, 
with $A(\tilde{u}_b) \in L^2\left(0,T; H^s(\Omega)\right)$.}
$$
 In fact, due to Theorem \ref{TraceTHM} 
this assumption could be sharpned by
$$
\text{ $u_b \in L^{\infty}(\Gamma_T) \cap L^2\left(0,T; H^{s-1/2}(\Gamma)\right)$,
$(\frac{1}{2}< s< 1)$.}
$$
\end{itemize}
Moreover, we denote for convenience  
$$
  a= \min \big\{ \ess \inf_{\Omega}(u_0),  \ess \inf_{\Gamma_T}(u_b)\big \} , \quad 
  b= \max \big\{ \ess \sup_{\Omega}(u_0),  \ess \sup_{\Gamma_T}(u_b) \big\},  
$$
and define 
\begin{equation}
\label{LFLA}
   L_f:= \ess\sup_{[a,b]} |f'|, \qquad
   L_A:= \ess\sup_{[a,b]} |A'|. 
\end{equation}

\medskip
The next lemma is important to establish the entropy formulation 
for the first equation in \eqref{FTPME}, where it is used the
Kru\v zkov's entropies.  
%
\begin{lemma}\label{lem:sgn}
Let $u$ be a smooth function in $\R^n$, and 
$A$ satisfying (H2).
Then, for each $k \in \R$ fixed,
and almost all $x \in \Omega$, 
$$
(-\Delta)^s_{\Omega} |A(u(x))-A(k)|^{\pm}= \sgn^{\pm}(u(x)-k)\,(-\Delta)^s_{\Omega} A(u(x))
-R_k(x), 
$$
where $R_k(\cdot)$ is a nonnegative function, given by 
$$
R_k(x):= C_{n,s} \int_{\Omega}\!\frac{(A(u(y))-A(k))(\sgn^{\pm}(u(x)-k)-\sgn^{\pm}(u(y)-k))}{|x-y|^{n+2s}} dy. 
$$
\end{lemma}
\begin{proof} First, 
observe that since $A(\cdot)$ is non-decreasing,
$$
\sgn^\pm\left(u(x)-k\right)\left(A(u(x))-A(k)\right)=\left|A(u(x))-A(k)\right|^\pm.
$$
Then, we have for almost all $x \in \Omega$
$$
\begin{aligned}
(-\Delta)^{s}_\Omega & |A(u(x))-A(k)|^{\pm}= 
(-\Delta)^{s}_\Omega \big( \sgn^{\pm} (u(x)-k) \, (A(u(x))-A(k)) \big)
\\[5pt]
&= \sgn^{\pm}(u(x)-k) \ (-\Delta)^{s}_\Omega \left( A(u(x)) -A(k)\right)
\\[5pt]
&+ \big(A(u(x))-A(k)\big) (-\Delta)^s_{\Omega}\, \sgn^{\pm}(u(x)-k)
\\[5pt]
&- C_{n,s} \int_\Omega\frac{\big(A(u(x))-A(u(y))\big) \, \big(\sgn^{\pm}(u(x)-k)-\sgn^{\pm}(u(y)-k)\big)}{|x-y|^{n+2s}} dy
\\[5pt]
&= \sgn^{\pm}(u(x)-k) \ (-\Delta)^{s}_\Omega A(u(x)) 
\\[5pt]
& +C_{n,s} \int_{\Omega}\!\frac{(A(u(y))-A(k))(\sgn^{\pm}(u(x)-k)-\sgn^{\pm}(u(y)-k))}{|x-y|^{n+2s}} dy,
\end{aligned}
$$
where we have used Lemma \ref{LFprodu}, together with $(-\Delta)^s_{\Omega}\,c=0$, when $c=const$. 
\end{proof}

\medskip
One recalls that, any smooth entropy pair 
$\mathbf{F}(u)= (\eta (u),q(u))$ for the first equation in \eqref{FTPME}, can be recovered by the family 
of Kru\v zkov's semi-entropies. Therefore, 
the following definition tell us in which sense a function $u \in L^\infty(Q_T)$ is a weak entropy solution of 
the initial-boundary value problem \eqref{FTPME}. 
\begin{definition}[Weak Entropy Solutions]
\label{DEFSOL}
Given $u_0 \in L^\infty(\Omega)$, 
$u_b \in L^\infty(\Gamma_T)$ satisfying $(H3)$, 
a function $ u\in L^{\infty}(Q_T)$ is said a weak entropy solution of  the problem \eqref{FTPME}, 
if it satisfies   
\begin{equation}
\label{eq:solfrac00}
    A(u)\in L^2\left(0,T; H^{s}(\Omega)\right), \quad A(u)-A(\tilde{u}_b) \in L^2(0,T;H^{s}_0(\Omega)),
\end{equation}
and  for each $k \in\R$ and any nonnegative test function $\psi \in C^{\infty}_c\left((-\infty,T)\times\R^n\right)$
such that, $\sgn^{\pm}\left(A(u_b) - A(k)\right) \, \psi= 0$ a.e. on $\Gamma_T$
\begin{equation}
\label{eq:solfrac}
\begin{aligned}
    \iint_{Q_T}& \Big\{ |u(t,x)-k|^{\pm} \, \partial_t\psi + \sgn^\pm(u(t,x)-k) (f(u(t,x)) - f(k)) \cdot \nabla\psi \Big\} \, dx dt 
    \\[5pt]
    &-\frac{C_{n,s}}{2}\int_0^T\!\!\left[|A(u(t,\cdot))-A(k)|^{\pm}\, , \,\psi(t,\cdot)\right]_{W^{s,2}(\Omega)} \, dt
  \\[7pt]  
    &+ L_f \int_{\Gamma_T}  |u_b(r) - k|^{\pm} \psi(r) \ d\mathcal{H}^{n}(r)
    +\int_\Omega  |u_0(x)-k|^{\pm} \psi(0) \, dx \geq 0,
\end{aligned}
\end{equation}
where $[\cdot,\cdot]_{W^{s,2}(\Omega)}$ is the Gagliardo semi-inner product, see \eqref{G_semi_inner_prod}.
\end{definition}

\begin{remark}
One remarks that, Definition \ref{DEFSOL} above is similar to Definition 4 in \cite{Carrillo}, 
and nearly to Definition 3.1 in \cite{Vovelle}, (see also Definition 2.1 in \cite{KobayasiOhwa}). 
Indeed, we recall that Carrillo considered homogeneous boundary 
condition, which together with the spaces considered by him for the test function, it is equivalent to 
$$\sgn^{\pm}\left(A(u_b) - A(k)\right) \, \psi= 0 \quad \text{a.e. on $\Gamma_T$}$$
considered here. 
\end{remark}

Henceforth, we extend for convenience the weak entropy solution $u \in L^\infty(Q_T)$
to $u \in L^\infty(\R \times \Omega)$ by setting  
$$
   u(t,x) \equiv 0, \quad \text{for all $(t,x) \in (\R \times \Omega) \setminus Q_T$}. 
$$

In the following we state the main theorem of the paper. 
\begin{theorem}[Main Theorem] 
\label{MainThm}
Let $f$ be the flux-function satisfying $(H1)$, and $A$ satisfying $(H2)$. 
Let $u_0 \in L^\infty(\Omega)$ be the initial data, 
$u_b \in L^\infty(\Gamma_T)$ be the boundary condition satisfying $(H3)$, 
and $\frac{1}{2}< s <1$. 
Then, there exists a function $u\in L^{\infty}(Q_T)$  
which is a weak entropy solution of the problem \eqref{FTPME}. 
Under the above conditions, 
if $u$ and $v$ are two weak entropy solutions of the problem \eqref{FTPME},
(obtained via parabolic perturbation), 
with respectively initial data $u_0$, $v_0$ and boundary condition $u_b$, $v_b$, such that
$A(\tilde{u}_b), A(\tilde{v}_b) \in L^2\left(0,T; H^{2s}(\Omega)\right)$, 
then for each $t \in (0,T)$,
\begin{equation}
\label{L1-contration}
\begin{aligned}
\int_{\Omega}|u(t,x)-v(t,x)| \, dx
&\leq\,\int_{\Omega}|u_0(x)-v_0(x)| \, dx
\\[5pt] 
&+ L_f \int_0^t\!\!\int_{\Gamma}|u_b(\tau,r)-v_b(\tau,r)|drd\tau
\\[5pt]
&-\mathcal{N}_\sigma\int_0^t\int_\Gamma\partial_\nu^\sigma|A(u_b)-A(v_b)|(\tau,r)drd\tau,  
\end{aligned}
\end{equation}
where $\sigma=2s-1$. Moreover, if $u_0 \leq v_0$ and $u_b \leq v_b$ almost everywhere, then 
$u(t,x) \leq v(t,x)$ for almost all $(t,x) \in Q_T$. 
\end{theorem}

\begin{remark}
\label{UsualContraction}
One observes that, if the boundary data is sufficiently regular, that is $C^1(\Gamma)$, then 
 \eqref{L1-contration} turns to the usual $L^1-$contraction, see Remark \ref{FracTrace}. 
\end{remark}

\section{On a parabolic perturbed problem}

The aim of this section is to introduce and study some properties of a 
parabolic perturbed equation associated to the problem \eqref{FTPME}, see \eqref{eq:regular}. More precisely, 
the existence of uniformly bounded solutions is considered for the viscous 
problem associated to \eqref{FTPME}. 
Hence an important inequality
is proved, which provides the existence
of weak entropy solutions of \eqref{FTPME} in the sense of Definition \ref{DEFSOL}.

\subsection{Parabolic approximation}

In order to show the existence of weak solutions for the 
initial-boundary value problem \eqref{FTPME}, we search
for $u_\ve$ be the solution,
for each parameter $\varepsilon>0$ fixed, of the following problem
\begin{equation}
\label{eq:regular}
\begin{cases}
   \partial_t u_\varepsilon+\varepsilon\,(-\Delta) u_\varepsilon+(-\Delta)^s_\Omega A_\ve(u_\varepsilon)= - \dive f(u_\varepsilon)  , 
   \quad &\mbox{ in }Q_T,
\\[5pt]
   u_\varepsilon= u_{0,\varepsilon},\quad &\mbox{ in }\Omega,
\\[5pt]
  u_\varepsilon= u_{b,\varepsilon}, \,\quad &\mbox{ on }\Gamma_T,
\end{cases}
\end{equation}
where  $f$, $A$ satisfy respectively conditions $(H1)$, $(H2)$, and $A_\ve(u):= A(u) + \ve u$.   Moreover, 
$u_{b,\varepsilon}$ and $u_{0,\varepsilon}$ are regularized boundary and initial data, respectively,
satisfying suitable compatibility conditions, such that 
\begin{equation}
\label{U0Ub}
\begin{aligned}
\text{$u_{0,\varepsilon} \to u_0 $ strongly in $L^{1}(\Omega)$
as $\varepsilon \to 0$,
and $\|u_{0,\varepsilon}\|_{L^\infty} \leq \|u_{0}\|_{L^\infty}$},
\\[5pt]
\text{$u_{b,\varepsilon} \to u_b$ strongly in $L^{1}(\Gamma_T)$
as $\varepsilon \to 0$,
and $\|u_{b,\varepsilon}\|_{L^\infty} \leq \|u_{b}\|_{L^\infty}$.}
\end{aligned}
\end{equation}
Following Lemmas 2,3 in \cite{Malek} and
due to Theorem \ref{Th_apen_Exis_2} from the Appendix,
we can show results of existence and uniqueness
for the initial-boundary value problem \eqref{eq:regular}. Moreover, making a 
regularization of the coefficients in $(\ref{eq:regular})_1$ and applying standard 
results of the regularity theory, we may assume that, for each $\ve> 0$ fixed, the
solution $u_\ve$ of \eqref{eq:regular} has enough regularity to 
perform the computations in Sections 4,5. Since $A_\ve$ converges uniformly to $A$ as $\ve \to 0^+$, 
for simplicity of notation, hereupon we continue to denote $A$ in \eqref{eq:regular} instead of $A_\ve$. 

\subsection{Perturbed problem estimates}

To begin, we establish a maximum principle for the solution $u_\ve$
of the initial-boundary value problem \eqref{eq:regular}. 

\begin{lemma}
\label{lema:Ref_PR}
Let  $u_\varepsilon$ be  the solution of \eqref{eq:regular}, hence it satisfies for each $t \in (0,T)$
\begin{eqnarray}
\int_{\Omega}\sgn^{-}(u_\varepsilon(t)-a_\ve)\,(-\Delta)_\Omega^s A(u_\varepsilon(t)) \ dx&\geq& 0,\nonumber
\\[5pt]
\int_{\Omega}\sgn^{+}(u_\varepsilon(t)-b_\ve)\,(-\Delta)_\Omega^s A(u_\varepsilon(t)) \ dx&\geq&0,\nonumber
\end{eqnarray}
where 
\begin{equation}
\label{abepsilon}
\begin{aligned}
a_\varepsilon:= \min\big( \inf_{\Omega}(u_{0,\varepsilon}) \ , \ \inf_{\Gamma_T}(u_{b,\varepsilon})\big),  \quad 
b_\ve:= \max\big( \sup_{\Omega}(u_{0,\varepsilon}) \ , \ \sup_{\Gamma_T}(u_{b,\varepsilon})\big). 
\end{aligned}
\end{equation}

\end{lemma}
\begin{proof}
Let us show the first inequality, the other one follows similarly. 
Applying Theorem \ref{Thintpart}, we obtain
$$
\begin{aligned}
   \int_{\Omega} \sgn^{-}&(u_\varepsilon(t,x) - a_\ve) \ (-\Delta)_\Omega^s A(u_\varepsilon(t,x)) \, dx
   \\[5pt]
   &= \frac{C_{n,s}}{2} \iint_{\Omega \times \Omega} \big(\sgn^{-}(u_\varepsilon(t,x) - a_\ve) - \sgn^{-}(u_\varepsilon(t,y) - a_\ve)\big)  
   \\[5pt]
    & \qquad \qquad \qquad \times \frac{A(u_\varepsilon(t,x)) - A(u_\varepsilon(t,y))}{|x-y|^{n+2s}}\,dx\,dy\,\geq 0,
\end{aligned} 
$$
where we have used that $\sgn^-(u_{b,\varepsilon} - a_\ve)= 0$. 
\end{proof}

\begin{proposition}[Maximum Principle]
\label{Ref_limita_u}
The solution $u_\varepsilon$ of \eqref{eq:regular} satisfies
\begin{equation}
\label{MAXPRINCIPLE1}
\sup_{Q_T} u_\varepsilon \leq b_\ve, \quad 
\inf_{Q_T} u_\varepsilon \geq a_\ve. 
\end{equation}
In particular,
\begin{equation}
\label{MAXPRINCIPLE2}
\sup_{Q_T} |u_\varepsilon| \leq \max  \big(\ess \sup_{\Omega} |u_{0}| \ , \ess \sup_{\Gamma_T} |u_{b}| \big).
\end{equation}
\end{proposition}
\begin{proof}
In order to prove the first inequality in \eqref{MAXPRINCIPLE1}, 
let $\varphi_\delta$ be defined by 
$$
\varphi_\delta(z):=\left\lbrace
\begin{aligned}
((z-b_\ve)^2+\delta^2)^{1/2}-\delta,&\quad, z\geq b_\ve,
\\[5pt]
0,&\quad,\, z\leq b_\ve.
\end{aligned}
\right.
$$
Then, we multiply $(\ref{eq:regular})_1$ by $\varphi^\prime_{\delta}(u_\ve)$, 
and using the properties of $\varphi_{\delta}$ we have
$$
\begin{aligned}
\int_{\Omega}& \varphi_{\delta}(u_\varepsilon(t)) \, dx 
- \int_{0}^{t} \int_{\Omega}\left({f}(u_\varepsilon)-{f}(b_\ve) \right) \cdot \nabla u_\varepsilon \ \varphi_{\delta}^{\prime \prime}(u_\varepsilon) \, dx d\tau 
\\[5pt]
&+\varepsilon \int_{0}^{t} \int_{\Omega}|\nabla u_\varepsilon|^{2} \varphi_{\delta}^{\prime \prime}(u_\varepsilon) \, dx d\tau
+\int_0^t\int_{\Omega}\varphi^{\prime}_\delta (u_\varepsilon) (-\Delta)^s_{\Omega} A(u_\varepsilon) \ dxd\tau= 0. 
\end{aligned}
$$
Now, we consider the estimate
$$
\begin{aligned}
-\left(  
{f}(u_\varepsilon)-{f}(b_\ve) \right) \cdot \nabla u_\varepsilon \ \varphi_{\delta}^{\prime \prime}(u_\varepsilon))&
+\varepsilon|\nabla u_\varepsilon|^{2} \varphi_{\delta}^{\prime \prime}(u_\varepsilon) 
\\[5pt]
& \geq\left\{-L_f|u_\varepsilon - b_\ve| |\nabla u_\varepsilon|+\varepsilon|\nabla u_\varepsilon|^{2}\right\} \varphi_{\delta}^{\prime \prime}(u_\varepsilon) 
\\[5pt]
& \geq-\frac{L_f^{2}}{4 \varepsilon}(u_\varepsilon - b_\ve)^{2} \varphi_{\delta}^{\prime \prime}(u_\varepsilon)
\geq-\frac{L_f^{2} \delta}{4 \varepsilon},
\end{aligned}
$$
and hence we obtain from the previous equation
$$
\int_{\Omega}\varphi_\delta(u_\varepsilon(t)) \, dx
+\int_0^t\int_{\Omega} \varphi^{\prime}_\delta (u_\varepsilon) (-\Delta)^s_{\Omega} A(u_\varepsilon) \, dx d\tau
\leq \frac{L_f^2\delta}{4\varepsilon}. 
$$
Therefore, passing to the limit as $\delta \rightarrow 0^+$ and using  Lemma \ref{lema:Ref_PR}, it follows that
$$
\int_{\Omega} \max (u_\varepsilon(t)-b_\ve, 0) \, dx \leq 0, 
$$
which is $(\ref{MAXPRINCIPLE1})_1$. Analogously, we obtain $(\ref{MAXPRINCIPLE1})_2$. 
Finally, the proof of $(\ref{MAXPRINCIPLE2})$ follows from $(\ref{MAXPRINCIPLE1})$ and \eqref{U0Ub}. 
\end{proof}

\medskip
The next theorem provides a suitable inequality satisfied by the solution 
$u_\ve$ of the problem \eqref{eq:regular}, which provides the existence
of weak entropy solutions of \eqref{FTPME} in the sense of Definition \ref{DEFSOL}.
\begin{theorem}
\label{EXISTHM}
Let $u_\ve$ be the unique solution of the problem \eqref{eq:regular}, 
and $\xi_\ve(x), L$ be defined by \eqref{XI4TECHOS}. 
Then, for each $k \in \R$ and for all nonnegative test function $\phi\in C^{\infty}_0((-\infty,T)\times\R^n)$, 
$u_\ve$ satisfies
$$
\begin{aligned}
&-\iint_{Q_T} \Big\{ \ \eta_k^\pm(u_{\varepsilon})\phi_t+ {q}_k^\pm(u_{\ve}) \cdot \nabla \phi- \ve \, \eta_k^\pm(u_{\ve})(-\Delta)\phi\Big\} \xi_{\varepsilon} \, dxdt
\\[5pt]
&+ \frac{C_{n,s}}{2} \iint_{Q_T} \!\! \xi_\ve(x) \!\! \int_{\Omega}   \big(\eta_{A(k)}^{\pm}(A(u_\ve(x))-\eta_{A(k)}^{\pm}(A(u_\ve(y)))\big) \frac{(\phi(x)-\phi(y))}{|x-y|^{n+2s}}\, dydx dt 
\\[5pt]
&\leq \int_{\Omega}\eta_k^\pm(u_{0,\varepsilon}) \, \phi(0) \xi_{\varepsilon} \, dx
+ (L_f + \varepsilon L) \int_{\Gamma_T} \eta_k^\pm(u_{b,\varepsilon}) \, \phi \, drdt
\\[5pt]
&+ 2 \varepsilon \iint_{Q_T} \eta_k^\pm(u_\ve) \nabla\phi \cdot \nabla\xi_\varepsilon \, dxdt
\\[5pt]
&-\frac{C_{n,s}}{2} \!\! \iint_{Q_T} \! \! \phi(t,x) \! \int_{\Omega}   \big(\eta_{A(k)}^{\pm}(A(u_\ve(x))-\eta_{A(k)}^{\pm}(A(u_\ve(y))\big)\frac{(\xi_\ve(x)-\xi_\ve(y))}{|x-y|^{n+2s}}\, dydx dt. 
\end{aligned}
$$
\end{theorem}
\begin{proof}
1. Let $(\eta_k^{\pm},q_k^{\pm})$ be 
the Kru\v zkov semi-entropy flux pair, and let $k \in \R$ 
be fixed, hence we denote $(\eta,{q}) \equiv (\eta_k^{\pm},{q}_k^{\pm})$.
To follow, we multiply equation $(\ref{eq:regular})_1$ by 
$\eta'(u_\varepsilon)$ and applying 
a standard procedure in conservation laws domain,
we obtain in distribution sense
\begin{equation}
\label{APPENTROP}
\begin{aligned}
\partial_t \eta (u_\varepsilon)  + \dive{q}(u_\varepsilon)
+ \varepsilon(-\Delta) \eta(u_\varepsilon)& + (-\Delta)_{\Omega}^s|A(u_\varepsilon)-A(k)|^{\pm}
\\[5pt]
&=-\varepsilon\eta''(u_\varepsilon)|\nabla u_\varepsilon|^2-R_k(x) \leq 0,
\end{aligned}
\end{equation}
where we have used Lemma \ref{lem:sgn}.

\medskip
2. Now, we multiply \eqref{APPENTROP}
by $\xi_{\varepsilon}(x) \, \phi(t,x)$, integrate over $Q_T$, and after integration by partes, we have
 (recall that $\xi_\ve \equiv 0$ on $\Gamma$)
\begin{equation}
   \iint_{Q_T} \partial_t \eta(u_\ve) \, \xi_\ve \phi \, dxdt
   =-\int_{\Omega} \eta(u_{0,\ve}) \,  \xi_\ve \phi(0) \, dx
   -\iint_{Q_T} \eta(u_\ve) \, \xi_\ve \, \partial_t\phi \, dxdt, 
\end{equation}
\begin{equation}
   \iint_{Q_T}\!\!\!\dive\,{q}(u_\ve) \, \xi_\ve \phi \, dxdt
= -\iint_{Q_T} \Big( \xi_\ve \, {q}(u_\ve) \cdot \nabla \phi 
   + \phi \, {q}(u_\ve) \cdot \nabla \xi_\ve\Big)\, dxdt, 
\end{equation}
\begin{equation}
\begin{aligned}
  \ve  \iint_{Q_T} &(-\Delta)\eta(u_\ve)\,\xi_\ve \phi \, dxdt
=\ve \iint_{Q_T}   \eta(u_\ve)\xi_\ve(-\Delta)\phi \, dx dt  
\\
&
+\ve \iint_{Q_T} \Big(\nabla\big(\eta(u_\ve)\phi\big) \cdot \nabla \xi_\ve  
 -2 \eta(u_\ve)\nabla\xi_\ve\cdot\nabla \phi \Big)\, dx dt. 
\end{aligned}
\end{equation}

\medskip
It remains to study the integral term, $\iint_{Q_T} \xi_\ve \phi \,  (-\Delta)^s_\Omega|A(u_\ve)-A(k)|^{\pm} \, dx dt$,
and to this end we observe that
\begin{equation}
\label{NEWTERM1}
\begin{aligned}
\iint_{Q_T}  \xi_\ve\phi& \, (-\Delta)^s_{\Omega}  \left(  \ |A(u_\ve)-A(k)|^{\pm}\right) \, dx dt
\\[5pt]
&=\frac{C_{n,s}}{2}\int_0^T\left[|A(u_\ve(t,\cdot))-A(k)|^{\pm}, \ \xi_\ve(\cdot)\phi(t,\cdot)\right]_{W^{s,2}(\Omega)}dt
\\
&\qquad - \mathcal{N}_\sigma \int_{\Gamma_T} \xi_\ve(r)\phi(t,r) \ \partial^\sigma_\nu \left( \ |A(u_\ve)-A(k)|^{\pm}\right) \, dr dt
\\
&=\frac{C_{n,s}}{2}\int_0^T\left[|A(u_\ve(t,\cdot))-A(k)|^{\pm}, \ \xi_\ve(\cdot)\phi(t,\cdot)\right]_{W^{s,2}(\Omega)}dt,
\end{aligned}
\end{equation}
where we have used Theorem \ref{Thintpart}, together with the fact that $\xi_\ve=0$ on $\Gamma$. 
Then, applying a simple algebraic manipulation, the last term in \eqref{NEWTERM1} can be written as
$$
\begin{aligned}
&\int_0^T\left[|A(u_\ve(t,\cdot))-A(k)|^{\pm}, \ \xi_\ve(\cdot)\phi(t,\cdot)\right]_{W^{s,2}(\Omega)}dt
\\[5pt]
&=\iint_{Q_T}\xi_\ve(x)\int_{\Omega}   \big(\eta_{A(k)}^{\pm}(A(u_\ve(x))-\eta_{A(k)}^{\pm}(A(u_\ve(y))\big)\frac{(\phi(x)-\phi(y))}{|x-y|^{n+2s}}\, dydx dt 
\\[5pt]
&+\iint_{Q_T}\phi(t,x)\int_{\Omega}   \big(\eta_{A(k)}^{\pm}(A(u_\ve(x))-\eta_{A(k)}^{\pm}(A(u_\ve(y))\big)\frac{(\xi_\ve(x)-\xi_\ve(y))}{|x-y|^{n+2s}}\, dydx dt.
\end{aligned}
$$

\medskip
3. Finally, we get together all the above terms to obtain
$$
\begin{aligned}
&-\iint_{Q_T} \Big\{ \ \eta(u_{\epsilon})\phi_t+ {q}(u_{\epsilon})\cdot\nabla \phi-\epsilon \eta(u_{\epsilon})(-\Delta)\phi \ \Big\} \xi_{\varepsilon} \, dxdt
\\[5pt]
&+\frac{C_{n,s}}{2} \iint_{Q_T} \!\! \xi_\ve(x)\!\! \int_{\Omega}  \big(\eta_{A(k)}^{\pm}(A(u_\ve(x))-\eta_{A(k)}^{\pm}(A(u_\ve(y))\big) \frac{(\phi(x)-\phi(y))}{|x-y|^{n+2s}}\,dy dx dt
\\[5pt]
&\leq \int_{\Omega}\eta(u_{0,\epsilon}) \, \phi(0) \xi_{\varepsilon} \, dx
+ \iint_{Q_T} \phi \, L_f  \, \eta(u_\varepsilon) |\nabla \xi_{\varepsilon}| \, dxdt
\\[5pt]
& -\ve \iint_{Q_T} \nabla(\phi \, \eta(u)) \cdot \nabla \xi_{\varepsilon} \, dxdt+ 2 \varepsilon \iint_{Q_T} \eta(u) \nabla\phi \cdot \nabla\xi_\varepsilon \, dxdt
\\[5pt]
&-\frac{C_{n,s}}{2} \!\! \iint_{Q_T} \!\! \phi(t,x) \!\! \int_{\Omega}   \big(\eta_{A(k)}^{\pm}(A(u_\ve(x))-\eta_{A(k)}^{\pm}(A(u_\ve(y))\big)\frac{(\xi_\ve(x)-\xi_\ve(y))}{|x-y|^{n+2s}}\, dydx dt,
\end{aligned}
$$
where we have used that $|{q}(\cdot)|\leq L_f \, \eta(\cdot)$. Therefore, applying inequality \eqref{preliminar}
with $\eta(u) \, \phi$ instead of $\beta$, we complete the proof.
\end{proof}

\section{Weak Entropy Solutions}

The main issue of this section is to show existence of weak entropy solutions to \eqref{FTPME}. 
To this end, we investigate more delicate properties of the solution $u_\ve$ to the problem
\eqref{eq:regular}. 
In particular, we proceed to show that the family $\{u_\ve\}$ of solutions to
 \eqref{eq:regular} is relatively compact in $L^1$. For that, we derive uniform estimates (with 
 respect to the parameter $\ve> 0$) on  
$$
   \|\partial_t u_\ve \|_{L^1(Q_T)} \quad \text{and} \quad   \|\nabla u_\ve \|_{L^1(Q_T)}. 
$$
We obtain such estimates taking the derivatives of 
the equation $(\ref{eq:regular})_1$, first with respect to the time 
variable, and then with respect to the spatial variable. 
One stress that, the 
Regional Fractional Laplacian term introduces many difficulties, 
which we were able to overcome. 

\medskip
Also, we prove a $L^1-$type contraction property for 
weak entropy solutions of \eqref{FTPME} obtained via parabolic 
perturbation. This is a weak selection principle which means that
the weak entropy solutions are stable w.r.t. the process of 
parabolic perturbation. 
In fact, the inequality obtained does not imply
necessarily stability, but in particular we obtain 
uniqueness of solution for equal data. Indeed, the stability condition is obtained for sufficiently smooth 
boundary data, see Remark \ref{UsualContraction}, 
therefore in this case,
we have well-posedness of the weak entropy solutions. 

\medskip
To begin, we drop for simplicity the subscript $\ve$ in 
$u_\ve$, $u_{0,\ve}$, $u_{b,\ve}$, and $a_\ve$, $b_\ve$ (see \eqref{abepsilon}), 
and let $u_{b}$ be the smooth extension (we keep the same notation) 
of $u_{b}$ onto $\overline{Q_T}$. 
Moreover, we fix the following notation
$$
\begin{aligned}
\|u_0\|_\Omega &\equiv \int_{\Omega} \big( |\nabla u_0| + |(-\Delta)u_0| 
+ |(-\Delta)^s_{\Omega} u_0| \big) dx, 
\\[5pt]
\|u_b\|_{\Gamma_T} &\equiv \int_{\Gamma_T} \big( |\partial_t u_b| + |\nabla u_b| + |\nabla \partial_t u_b| +  |(-\Delta)^s_\Omega u_b| + |D^2 u_b| \big) \, dr dt,
\\[5pt]
   \|u_b\|_{Q_T} &\equiv  \iint_{Q_T} \big(| \partial_t u_b| + |\nabla u_b| + |\Delta u_b| + |(-\Delta)^s_\Omega u_b| \big) \, dx dt. 
\end{aligned}
$$
 
We start with the uniform estimate with respect to the time derivative.    
\begin{proposition}
\label{TimeDerivativeEstimate}
Let $u$ be the solution of the problem \eqref{eq:regular} 
with initial data $u_0$ and boundary condition $u_b$ respectively. Then
\begin{equation}
\label{ESTIME}
    \sup_{t \in (0,T)} \int_{\Omega}\left |\partial_t u(t,x) \right|dx\leq C_1, 
\end{equation} 
where $0< C_1= C_1 \left(\|u_0\|_\Omega,\|u_b\|_{\Gamma_T}, \Omega, T, L_f, L_A\right)$ does not depend on $\varepsilon> 0$. 
\end{proposition}
\begin{proof} 1. First, denoting
$v= \partial_t u$
we obtain from equation $(\ref{eq:regular})_1$, 
\begin{equation}
\label{eq:estimativa_u_derivada_de_t}
\frac{\partial v}{\partial t}+\dive\left(f^{\prime}(u) v\right)-\varepsilon \Delta v+(-\Delta)^s_{\Omega}\left(A'(u)v\right)= 0.
\end{equation}
Then, multiplying \eqref{eq:estimativa_u_derivada_de_t} by $\varphi_{\delta}^{\prime}(v)$, 
where $\varphi_{\delta}(z) \equiv\left(z^{2}+\delta^{2}\right)^{1 / 2}$, and integrating over $(0, t) \times \Omega$, we have
\begin{equation}
\begin{aligned}
&\int_{\Omega} \varphi_{\delta}(v(t)) \, d x + \int_0^t \!\! \int_\Gamma \partial_t u_b \, \varphi_\delta^\prime(\partial_t u_b) \, f^\prime(u_b) \cdot \nu \, dr d\tau
\\[5pt]
& -\int_{0}^{t} \!\! \int_{\Omega} v \, \varphi_{\delta}^{\prime \prime}(v) \,f^{\prime}(u) \cdot \nabla v \, dx d\tau
- \ve \int_0^t \!\! \int_\Gamma \varphi_\delta^\prime(\partial_t u_b) \, \nabla(\partial_t u_b) \cdot \nu \, dr dt
\\[5pt]
&+ \varepsilon \int_{0}^{t} \!\! \int_{\Omega}|\nabla v|^{2} \varphi_{\delta}^{\prime \prime}(v) \, dx d\tau
+\frac{C_{n,s}}{2} \int_0^t \left[\varphi'_\delta(v),A'(u) v \right]_{W^{s,2}(\Omega)} \, d\tau 
\\[5pt]
&= \int_{\Omega} \varphi_{\delta}(v(0)) \, dx,
\end{aligned}\nonumber
\end{equation}
where we have used Remark \ref{FracTrace}, and 
$[\cdot,\cdot]_{W^{s,2}(\Omega)}$ is the Gagliardo semi-inner product.
Again, we consider the standard inequality 
$$
\begin{aligned}
-v \varphi_{\delta}^{\prime \prime}(v) \ f^{\prime}(u) \cdot \nabla v +\varepsilon|\nabla v|^{2} \varphi_{\delta}^{\prime \prime}(v) & \geq-\frac{1}{4 \varepsilon}\left|f^{\prime}(u)\right|^{2} v^{2} \varphi_{\delta}^{\prime \prime}(v) \\
& \geq-\frac{\delta \, L_f^{2}}{4 \varepsilon}
\end{aligned}
$$
and thus, it follows from the above equation 
\begin{equation}
\label{eq2:Pro_estimativa}
\begin{aligned} 
&\int_{\Omega} \varphi_{\delta}(v(t)) \, dx  + \int_0^t \!\! \int_\Gamma \partial_t u_b \, \varphi_\delta^\prime(\partial_t u_b) \, f^\prime(u_b) \cdot \nu \, dr d\tau
\\[5pt]
& - \ve \int_0^t \!\! \int_\Gamma \varphi_\delta^\prime(\partial_t u_b) \, \nabla(\partial_t u_b) \cdot \nu \, dr dt
+\frac{C_{n,s}}{2}\int_0^t \left[\varphi'_\delta(v),A'(u) v \right]_{W^{s,2}(\Omega)} \, d\tau 
\\[5pt]
&\leq \int_{\Omega} \varphi_{\delta}(v(0)) \, d x+\frac{\delta \, L_f^{2}\, T}{4 \varepsilon} \, |\Omega|. 
\end{aligned}
\end{equation}
Therefore, passing to the limit as $\delta \rightarrow 0^+$ in \eqref{eq2:Pro_estimativa}, we have
$$
\begin{aligned}
& \int_{\Omega}|v(t)|\, d x+\frac{C_{n,s}}{2}\int_0^t \left[\sgn(v(\tau)),A'(u(\tau))v(\tau) \right]_{W^{s,2}(\Omega)}d\tau
\\[5pt]
&\leq \int_{\Omega}|v(0)| \, d x+\int_{0}^{t} \!\! \int_{\Gamma} \big( L_f |\partial_t u_b| + |\nabla \partial_t u_b| \big) \, dr d\tau. 
\end{aligned}
$$

2. Now, we consider the following 

\smallskip
 \underline {Claim:} For each $t \in (0,T)$, 
 \begin{equation}
 \label{Lem:produto_inteno_positivo}
     \left[\sgn(v(t)),A'(u(t))v(t) \right]_{W^{s,2}(\Omega)} \geq 0. 
\end{equation} 

Proof of Claim: Denoting $u(t,\cdot)\equiv u(\cdot)$, and $v(t,\cdot)\equiv v(\cdot)$, it follows by definition that
$$
\begin{aligned}
&\left[\sgn(v(t)),A'(u(t)) \, v(t)\right]_{W^{s,2}(\Omega)}
\\[5pt]
&=\iint_{\Omega \times \Omega}\!\!\! \!\! \frac{(A'(u(x)) v(x) - A'(u(y)) v(y))(\sgn(v(x)) - \sgn(v(y))}{|x-y|^{n+2s}}dxdy \geq 0.
\end{aligned}
$$
Indeed, clearly if $v(x)$ and $v(y)$ have the same sign, then 
the above integral is zero. On the other hand, if they have different signs 
jointly with $A^\prime(\cdot) \geq 0$, then the result follows.   

To follow, we observe that 
$$
\begin{aligned}
\int_{\Omega}|v(0)| \, dx&= \int_{\Omega}\left|\frac{\partial u}{\partial t}(0)\right| \, dx
=\int_{\Omega}\left|-\dive f \left(u_{0}\right)+\varepsilon \Delta u_{0}-(-\Delta)^s_{\Omega}A(u_0)\right| \, dx
\\[5pt]
&\leq L_f\int_{\Omega}|\nabla u_0| dx + \int_{\Omega}|(-\Delta)u_0| dx 
+ L_A \int_{\Omega}|(-\Delta)^s_{\Omega} u_0| dx. 
\end{aligned}
$$

3. Finally, from items $(1)$ and $(2)$ we obtain for each $t \in (0,T)$ 
$$
\begin{aligned}
\int_{\Omega} \left|\frac{\partial u(t)}{\partial t}\right| d x &\leq L_f\int_{\Omega}|\nabla u_0| dx + \int_{\Omega}|(-\Delta)u_0| dx 
+ L_A  \int_{\Omega}|(-\Delta)^s_{\Omega} u_0| dx 
\\[5pt]
& +\int_{0}^{t} \!\! \int_{\Gamma} \big( L_f |\partial_t u_b| + |\nabla \partial_t u_b| \big) \, dr dt, 
\end{aligned}
$$
from which the result follows.
\end{proof}

The main issue now is to establish the uniform estimate with respect to the spatial derivative.
Let us begin considering three important lemmas. The former one is
\begin{lemma} 
\label{LEMMADERIVX}
Let $w$ be a smooth scalar function in $\R^n$. Then, $w$ satisfies for each $i \in (1,\ldots, n)$, 
$$
\begin{aligned}
\frac{\partial}{\partial x_i}\left((-\Delta)^s_{\Omega}w\right)(x)=(-\Delta)^s_{\Omega}\left(\frac{\partial w}{\partial x_i}\right)(x)
+C_{n,s}\int_{\Gamma}\frac{w(r)-w(x)}{|r-x|^{n+2s}}\nu_i(r) \ dr.
\end{aligned}
$$
\end{lemma}
\begin{proof} Let $\epsilon> 0$ be fixed, also $i \in (1,\ldots, n)$.
Hence we have by definition 
$$
\begin{aligned}
&\frac{\partial}{\partial x_i}\left((-\Delta)^s_{\Omega,\epsilon}w\right)(x)
=\,C_{n,s} \ \frac{\partial}{\partial x_i}\left(\int_{\Omega\setminus B_\epsilon(x)}\frac{w(x)-w(y)}{|x-y|^{n+2s}}dy\right)
\\[7pt]
&=C_{n,s}\int_{|r-x|=\epsilon}\frac{w(x)-w(r)}{|x-r|^{n+2s}}\nu_i(r) \ dr+\,C_{n,s}\int_{\Omega\setminus B_\epsilon(x)}\frac{\partial}{\partial x_i}\left(\frac{w(x)-w(y)}{|x-y|^{n+2s}}\right) \, dy
\\[7pt]
&=C_{n,s}\int_{|r-x|=\epsilon}\frac{w(x)-w(r)}{|x-r|^{n+2s}}\nu_i(r) \ dr+\,C_{n,s}\int_{\Omega\setminus B_\epsilon(x)}\frac{\partial_{x_i}w(x)}{|x-y|^{n+2s}} \, dy
\\[7pt]
&-C_{n,s}\int_{\Omega\setminus B_\epsilon(x)}\big(w(x)-w(y)\big)\frac{\partial_{x_i}(|x-y|^{n+2s})}{|x-y|^{2(n+2s)}} \, dy. 
\\[7pt]
\end{aligned}
$$
One observes that, $\partial_{x_i}(|x-y|^{n+2s})=-\partial_{y_i}(|x-y|^{n+2s})$, 
then we have 
$$
\begin{aligned}
&\frac{\partial}{\partial x_i}\left((-\Delta)^s_{\Omega,\epsilon}w\right)(x)
\\
&=C_{n,s}\int_{|r-x|=\epsilon}\frac{w(x)-w(r)}{|x-r|^{n+2s}}\nu_i(r) \ dr+\,C_{n,s}\int_{\Omega\setminus B_\epsilon(x)}\frac{\partial_{x_i}w(x)}{|x-y|^{n+2s}} \, dy
\\[5pt]
&+C_{n,s}\int_{\Omega\setminus B_\epsilon(x)}\big(w(x)-w(y)\big)\frac{\partial_{y_i}(|x-y|^{n+2s})}{|x-y|^{2(n+2s)}}dy
\\[5pt]
&=C_{n,s}\int_{|r-x|=\epsilon}\frac{w(x)-w(r)}{|x-r|^{n+2s}}\nu_i(r) \ dr+\,C_{n,s}\int_{\Omega\setminus B_\epsilon(x)}\frac{\partial_{x_i}w(x)}{|x-y|^{n+2s}} \, dy
\\[5pt]
&-C_{n,s}\int_{\Omega\setminus B_\epsilon(x)}\big(w(x)-w(y)\big)\partial_{y_i}\left(\big(|x-y|^{n+2s}\big)^{-1}\right) dy.
\\[5pt]
\end{aligned}
$$
Therefore, we make a partial integration in the third term to get
$$
\begin{aligned}
\frac{\partial}{\partial x_i}&\left((-\Delta)^s_{\Omega,\epsilon}w\right)(x)
\\
&=\,C_{n,s}\int_{\Omega\setminus B_\epsilon(x)}\frac{\partial_{x_i}w(x)-\partial_{y_i}w(y)}{|x-y|^{n+2s}}dy
+C_{n,s}\int_{\Gamma}\frac{w(r)-w(x)}{|x-r|^{n+2s}}\nu_i(r) \, dr. 
\end{aligned}
$$
To finish, we pass to the limit as $\epsilon\to0^+$, (observe that this limit is uniform),
thus the result follows.
\end{proof}

\begin{lemma} 
\label{lem:estimati_positivo}
Let $\mathbf{w}= (w_1, w_2, \ldots, w_n)$ be a smooth vector field in $\R^n$, such that $|\mathbf{w}|> 0$.
Then, 
$$
\sum_{i=1}^n\left[\frac{w_i}{|\mathbf{w}|}\, , \, A'(\cdot)w_i\right]_{W^{s,2}(\Omega)} \geq 0. 
$$
\end{lemma}
\begin{proof}
By definition, we need to show that 
$$
\iint_{\Omega\times\Omega}
\Big(\frac{w_i(x)}{|\mathbf{w}(x)|}-\frac{w_i(y)}{|\mathbf{w}(y)|}\Big)
\big( A'(u(x)) w_i(x)-A'(u(y)) w_i(y)\big)\frac{1}{|x-y|^{n+2s}}\geq 0,
$$
where $u$ is a smooth function in $\R^n$ and we are using the summation convention. 
For any $x,y \in \Omega$, we have 
$$
\begin{aligned}
&\left(\frac{w_i(x)}{|\mathbf{w}(x)|}-\frac{w_i(y)}{|\mathbf{w}(y)|}\right)\left( A'(u(x)) w_i(x) - A'(u(y)) w_i(y)\right)
\\[5pt]
&=A'(u(x))|\mathbf{w}(x)|+A'(u(y))|\mathbf{w}(y)|
- w_i(x)w_i(y)\left(\frac{A'(u(x))}{|\mathbf{w}(y)|}+\frac{A'(u(y))}{|\mathbf{w}(x)|}\right)
\\[5pt]
&=\left(|\mathbf{w}(y)||\mathbf{w}(x)|
- w_i(x)w_i(y)\right)\left(\frac{A'(u(x))}{|\mathbf{w}(y)|}+\frac{A'(u(y))}{|\mathbf{w}(x)|}\right).
\end{aligned}
$$
Applying the Cauchy-Swartz inequality and due to $A'(\cdot) \geq 0$, the result follows.
\end{proof}

Now for $\rho>0$, we consider the following function
\begin{equation}
\label{eq:beta_rho}
\beta_\rho(x):=\gamma\left(\frac{h(x)}{\rho}\right),\quad x\in\R^n,
\end{equation}
where $h(x)$ is defined in the Appendix and $\gamma \in C^{\infty}(\mathbb{R})$ is a fixed non-negative function, such that
$$
\gamma(0)=0, \quad \gamma(\theta)= 1 \quad \text { for } \theta \geq 1.
$$
Also, for $\delta > 0$ (small enough), we recall the definition of $\Omega_\delta$, see \eqref{Omega_delta}, and define the set 
$\Omega_{\delta}^c$, respectively, 
$$
\begin{aligned}
\Omega_{\delta}&:=\left\lbrace x\in\Omega:\, 0< d(x)< \delta\right\rbrace,
\\[5pt]
\Omega_{\delta}^c&:= \Omega\setminus\Omega_\delta= \left\lbrace x \in \Omega : \, \delta \leq d(x) \right\rbrace.
\end{aligned}
$$
Then, we have the following 
\begin{lemma} 
\label{lem:Lap_lim_A_beta} 
Let $\Psi$ be a smooth function in $\R^n$, and $\beta_\rho$ as defined in \eqref{eq:beta_rho}. Then, 
$$
\begin{aligned}
\lim_{\rho\to0^+}\int_{\Omega}\Psi(x)(-\Delta)^s_{\Omega}\beta_{\rho}(x)dx&=C_{n,s}\int_{\Omega}\int_{\Gamma}\frac{\Psi(x)-\Psi(r)}{|x-r|^{n+2s}} \, dx dr. 
\end{aligned}
$$
\end{lemma}
\begin{proof}
1. First, applying the Gauss-Green Theorem we have
$$
\begin{aligned}
\int_{\Omega}\Psi(x)(-\Delta)^s_{\Omega}\beta_{\rho}(x)dx&=\frac{C_{n,s}}{2}
\iint_{\Omega \times \Omega} \frac{\left(\Psi(x)-\Psi(y)\right)\left(\beta_{\rho}(x)-\beta_{\rho}(y)\right)}{|x-y|^{n+2s}}dxdy
\\[5pt]
&=
\frac{C_{n,s}}{2}
\iint_{\Omega_\delta \times \Omega_\delta} \frac{\left(\Psi(x)-\Psi(y)\right)\left(\beta_{\rho}(x)-\beta_{\rho}(y)\right)}{|x-y|^{n+2s}}dxdy
\\[5pt]
&+ C_{n,s}
\iint_{\Omega_\delta \times \Omega_\delta^c}\frac{\left(\Psi(x)-\Psi(y)\right)\left(\beta_{\rho}(x)-\beta_{\rho}(y)\right)}{|x-y|^{n+2s}}dxdy
\\[5pt]
&=: I_1 + I_2
\end{aligned}
$$
with the obvious notation, where we have used that, for each $(x,y) \in \Omega_{\delta}^c \times \Omega_{\delta}^c$ we have  
$h(x)= h(y)= \delta$. 

\medskip
2. Now, let us study $I_2$. 
We consider that $\rho>0$ is small enough, such that $\rho<\delta$, that is,
$1< \delta / \rho$. Consequently, $\beta_{\rho}(x)=1$ for all $x\in\Omega_{\delta}^c$. 
To follow, since $\Gamma$ is a $C^2-$boundary, each $y\in\Omega_\delta$ 
 has unique projection $r=\mathbf{r}(y)$ on the boundary, and the Jacobian of the change of variables,
$$\Omega_{\delta} \ni y \leftrightarrow ({r},\tau) \in \Gamma \times (0,\delta),$$
is $1 + O(\delta )$, where $\tau= \dist(y,\Gamma)$. Therefore, we can write 
$$
\begin{aligned}
     I_2&= C_{n,s}\int_{\Omega_{\delta}} \! \int_{\Omega_{\delta}^c}\frac{\left(\Psi(x)-\Psi(y)\right)\left(1-\beta_{\rho}(y)\right)}{|x-y|^{n+2s}}dxdy
\\[5pt]
&=C_{n,s} \ \rho\int_0^{\frac{\delta}{\rho}} \! \! \int_{\Gamma}\int_{\Omega_\delta^c}\frac{\big(\Psi(x)-\Psi(r-\tau\rho\nu)\big)(1-\gamma(\tau))}{|x-\left(r-\tau\rho\nu\right)|^{n+2s}}dydrd\tau+O(\delta).
\end{aligned}
$$
Then passing to the limit as $\rho\to0^+$, and after that making $\delta \to 0^+$, we get that the $I_2$ term converges to
$$
\begin{aligned}
C_{n,s} \ \int_{\Gamma}\int_{\Omega}\frac{\Psi(x)-\Psi(r)}{|x-r|^{n+2s}}\,dxdr. 
\end{aligned}
$$

\medskip
3. Finally, we study the $I_1$ term. We make the similar change of variables 
as done in item 2, then we have
$$
\begin{aligned}
I_1&= \frac{C_{n,s}}{2}
\iint_{\Omega_\delta \times \Omega_\delta} \frac{\left(\Psi(x)-\Psi(y)\right)\left(\beta_{\rho}(x)-\beta_{\rho}(y)\right)}{|x-y|^{n+2s}}dxdy
\\[5pt]
& =\frac{C_{n,s}}{2}\rho^2 \!\! \int_0^\frac{\delta}{\rho}\!\! \int_0^\frac{\delta}{\rho}
 \iint_{\Gamma\times\Gamma} \frac{\big(\Psi(r_x - \tau_x \rho \nu(r_x)) - \Psi(r_y- \tau_y \rho\nu(r_y)\big)}{|r_x - \tau_x \rho \nu(r_x)-(r_y-\tau_y\rho\nu(r_y))|^{n+2s}}
 \\[5pt]
  & \hspace{130pt}  \times \big(\gamma(\tau_x)-\gamma(\tau_y)\big) \, dr_x dr_y d\tau_x d\tau_y + O(\delta), 
\end{aligned}
$$
with the obvious notation. Then, passing to the limit as $\rho\to0^+$, 
and making $\delta \to 0^+$, it follows that $I_1$ converges to zero, which concludes the proof.
\end{proof}

\medskip
Now, we are able to establish the uniform estimate with respect to the spatial derivatives.    
\begin{theorem}
\label{pro:estimativa_w}
Let $u$ be the solution of the problem \eqref{eq:regular} 
with initial data $u_0$ and boundary condition $u_b$ respectively. Then, for each $t \in (0,T)$
$$
\begin{aligned}
(&i) \int_{\Omega}|\varpi(t)| d x + \! \varepsilon \!\! \int_{0}^{t}\! \! \int_{\Gamma}|\nabla \varpi \cdot \nu| d r d \tau 
\\[5pt]
& +C_{n,s} \! \! \int_0^t \! \! \int_{\Omega} \! \int_{\Gamma}\frac{|A(u(\tau,x))-A(u_b(\tau,x))|}{|x-r|^{n+2s}}drdxd\tau
\leq \int_\Omega \varpi(0) \ dx + \!\! \int_{0}^{t}\!\! \int_{\Omega}|g| d x d \tau,
\end{aligned}
$$
where 
$\varpi(t,x):= u(t,x) - u_b(t,x)$, and 
$$
    g:= \partial_t u_{b}+\dive\, f\left(u_{b}\right)-\varepsilon \Delta u_{b}+(-\Delta)^s_{\Omega}A(u_b).
$$    
Moreover, there exists $0< C_2= C_2 \left(\|u_0\|_\Omega,\|u_b\|_{\Gamma_T}, \|u_b\|_{Q_T}, \Omega,T,L_f, L_A, a, b\right)$, 
which does not depend on $\varepsilon> 0$, such that
\begin{equation}
\label{ESTX}
  (ii) \quad \sup_{t\in(0,T)} \int_{\Omega} \left|\nabla u(t,x)\right|dx \,\leq\, C_2. \hspace{156pt} 
\end{equation}
\end{theorem}
\begin{proof}
1. First, let us show item $(i)$. 
From equation $(\ref{eq:regular})_1$, we obtain
\begin{equation}
\label{pro:estimativa_w}
\frac{\partial \varpi}{\partial t}+\dive\left(f(u)-f\left(u_{b}\right)\right)-\varepsilon \Delta \varpi + (-\Delta)^s_{\Omega}\big(A(u)-A(u_b)\big)=- g.
\end{equation}
Then, we multiply\eqref{pro:estimativa_w} by $\varphi_{\delta}^{\prime}(\varpi) \beta$, with $\beta \in C^\infty_c(\mathbb{R}^{d})$, $\beta \geq 0$, and
$$
\varphi_{\delta}(z)= \left(z^{2} + \delta^{2}\right)^{1 / 2} - \delta,
$$
integrate over $(0, t) \times \Omega$, and after partial integration we have for each $t \in (0,T)$, 
$$
\begin{aligned}
&\int_{\Omega} \varphi_{\delta}(\varpi(t)) \beta \, d x - \!\!  \int_{\Omega} \varphi_{\delta}(\varpi(0)) \beta \, d x 
- \!\! \int_{0}^{t}\!\! \int_{\Omega} \varphi_{\delta}^{\prime}(\varpi(\tau)) \left(f(u)-f(u_b)\right) \cdot \nabla \beta \, dx d\tau 
\\[5pt]
&-\int_{0}^{t} \!\! \int_{\Omega} \varphi_{\delta}^{\prime \prime}(\varpi(\tau)) \beta \left(f(u)-f(u_b)\right) \cdot \nabla \varpi(\tau) \, dx d\tau 
\\[5pt]
&+\varepsilon \int_{0}^{t} \!\! \int_{\Omega} |\nabla \varpi(\tau)|^{2} \varphi_{\delta}^{\prime \prime}(\varpi(\tau)) \beta \, d x d \tau
-\varepsilon \int_{0}^{t} \!\! \int_{\Omega} \varphi_{\delta}(\varpi(\tau)) \Delta \beta(x) d x d \tau 
\\[5pt]
&+\int_0^t \beta \, \varphi_\delta^{\prime}(\varpi(t)) \, (-\Delta)^s_{\Omega}\big(A(u)-A(u_b)\big)
=-\int_{0}^{t} \!\! \int_{\Omega} \varphi_{\delta}^{\prime}(\varpi(\tau)) g(\tau) \beta(x) d x d \tau,
\end{aligned}
$$
where we have used that
$\varphi_{\delta}^{\prime}(\varpi) \equiv \varphi_{\delta}(\varpi) \equiv \nabla \varphi_{\delta}(\varpi) \cdot \nu= 0$ on $\Gamma$. 
Then, passing to the limit as $\delta \rightarrow 0^+$, it follows that 
\begin{equation}
\label{XESTIM}
\begin{aligned}
&\int_{\Omega}|\varpi(t)|\beta \, d x - \int_{\Omega}|\varpi(0)|\beta \, d x 
-\int_{0}^{t} \!\! \int_{\Omega} q(u,u_b) \cdot \nabla \beta \, d x d \tau 
\\[5pt]
&-\varepsilon \int_{0}^{t}\!\!  \int_{\Omega} |\varpi| \, \Delta \beta \, d x d \tau
+ \int_0^t \!\! \int_{\Omega} |A(u)-A(u_b)| \, (-\Delta)^s_{\Omega} \beta \, dx d\tau
\\[5pt]
&\leq-\int_{0}^{t} \!\! \int_{\Omega} \sgn(\varpi(\tau)) \ g(\tau) \, \beta \, d x d \tau, 
\end{aligned}
\end{equation}
where $q(u,u_b)=\sgn(u-u_b)\big(f(u)-f(u_b)\big)$.

\medskip
2. Now, we choose in \eqref{XESTIM}, $\beta= \beta_\rho$ (see \eqref{eq:beta_rho}), and 
performing as item (3) in the proof of Lemma \ref{lem:Lap_lim_A_beta}, we observe that 
$$
\begin{aligned}
&\lim _{\rho \rightarrow 0+} \int_{0}^{t}\!\!  \int_{\Omega} q\left(u, u_{b}\right) \nabla \beta_{\rho} \, d x d \tau 
\\[5pt]
&=\lim _{\rho \rightarrow 0+} \int_{0}^{\frac{\delta}{\rho}} \gamma^{\prime}(\sigma) \int_{0}^{t} \! 
\int_{\Gamma} q(u, u_{b})(\tau, r - \rho \sigma\nu(r)) \cdot \nu(r) d r d \tau d \sigma + O(\delta)
\\[5pt]
&= O(\delta)= 0, 
\end{aligned}
$$
where we have used that $\beta_{\rho}= 0$ on $\Gamma$, and sent $\delta \to 0^+$. Analogously, we have 
$$
\begin{aligned}
\lim _{\rho \rightarrow 0+} \int_{0}^{t} \!\! \int_{\Omega}|\varpi| & \Delta \beta_{\rho} \, d x d \tau
=-\int_{0}^{t} \int_{\Gamma}|\nabla \varpi \cdot \nu| d r d \tau. 
\end{aligned}
$$

3. Due to Lemma \ref{lem:Lap_lim_A_beta} and taking into account that $u=u_b$ on $\Gamma$, it follows 
for each $t \in (0,T)$
$$
\lim_{\rho\to0^+}\int_{\Omega}|A(u)-A(u_b)|(t)(-\Delta)^s_{\Omega}\beta_{\rho}dx=
C_{n,s} \!\! \int_{\Omega} \! \int_{\Gamma}\frac{|A(u)-A(u_b)|(t)}{|x-r|^{n+2s}} \, dxdr,
$$
from which follows the proof of item $(i)$. 

\medskip 
4. Now, we prove item $(ii)$. Let us denote $\mathbf{w}= (w_{1}, w_2, \ldots, w_{n})$, where  
for $(t,x) \in Q_T$, 
$$
    w_{i}(t,x)= \frac{\partial u(t,x)}{\partial x_{i}}, \qquad  (i=1, \ldots, n). 
$$
Therefore, we obtain from equation $(\ref{eq:regular})_1$, for each $i \in (1, \ldots,n)$ 
\begin{equation}
\label{pro:equa_w}
\begin{aligned}
\frac{\partial w_{i}(t,x)}{\partial t}+\dive\left(f^{\prime}(u) w_{i}\right)&-\varepsilon \Delta w_{i}+(-\Delta)^s_{\Omega}(A'(u)w_i)
\\[5pt]
&+C_{n,s}\int_{\Gamma}\frac{A(u_b(t,r))-A(u(t,x))}{|r-x|^{n+2s}}\nu_i(r)dr= 0,
\end{aligned}
\end{equation}
where we have used Lemma \ref{LEMMADERIVX}.
Then, we multiply the equation \eqref{pro:equa_w} by $\frac{\partial}{\partial \xi_{i}} \phi_{\delta}(\mathbf{w})$, where
$\phi_{\delta}(\xi)=\left(|\xi|^{2}+\delta^{2}\right)^{1 / 2}$, and integrate over $(0, t) \times \Omega$. Let
us study each term separately, where we use the summation convention, we have
$$
\begin{aligned}
\int_{0}^{t} \!\! \int_{\Omega} \frac{\partial w_{i}}{\partial t} \, \frac{\partial \phi_{\delta}}{\partial \xi_{i}}(\mathbf{w}) \, d x d \tau 
\,&=\int_{\Omega} \phi_{\delta}(\mathbf{w}(t))\, d x - \int_{\Omega} \phi_{\delta}(\mathbf{w}(0)) \, d x, 
\\[5pt]
-\varepsilon \int_{0}^{t} \!\! \int_{\Omega} \Delta w_{i}  \frac{\partial \phi_{\delta}}{\partial \xi_{i}}(\mathbf{w}) \,d x d \tau 
\,& = \varepsilon \int_{0}^{t} \!\! \int_{\Omega} \frac{\partial w_i}{\partial x_j} \frac{\partial^2 \phi_\delta}{\partial \xi_i \partial \xi_k}(\mathbf{w}) \frac{\partial w_k}{\partial x_j}  \, d x d \tau 
\\[5pt]
&-\varepsilon \int_{0}^{t}  \! \! \int_{\Gamma}  \frac{\partial \phi_{\delta}}{\partial \xi_{i}}(\mathbf{w}_b) \ \nabla w_{b i} \cdot \nu(r) \, d r d \tau,
\\[5pt]
\int_{0}^{t} \!\! \int_{\Omega} \dive\left(f^{\prime}(u) w_{i}\right) \frac{\partial \phi_{\delta}}{\partial \xi_{i}}(\mathbf{w}) \, d x d \tau 
&=-\int_{0}^{t} \!\! \int_{\Omega}w_i \, f'_j(u) \frac{\partial^2 \phi_\delta}{\partial \xi_i \partial \xi_k}(\mathbf{w}) \frac{\partial w_k}{\partial x_j} \, d x d \tau 
\\[5pt]
&+\int_{0}^{t} \int_{\Gamma}  w_{bi} \frac{\partial \phi_{\delta}}{\partial \xi_{i}}(\mathbf{w}_b) f'(u_b) \cdot \nu(r) \, d r d \tau,
\end{aligned}
$$
with obvious notation.  Moreover, we have for each $t \in (0,T)$
$$
\begin{aligned}
\int_{\Omega} \frac{\partial \phi_{\delta}}{\partial \xi_{i}}(\mathbf{w}(t)) &\, (-\Delta)^s_{\Omega}(A'(u(t))w_i(t)) \, dx
\\
&= \frac{C_{n,s}}{2} \left[\frac{\partial \phi_{\delta}}{\partial \xi_{i}}(\mathbf{w}(t)),A'(u(t)) w_i(t) \right]_{W^{s,2}(\Omega)},
\end{aligned}
$$
where we have used, again, Remark \ref{FracTrace}. 
We also consider the following standard estimate
$$
\begin{aligned}
& \ve \,  \frac{\partial w_i}{\partial x_j}  \frac{\partial^2 \phi_\delta}{\partial \xi_i \partial \xi_k}(\mathbf{w}) \frac{\partial w_k}{\partial x_j}  
- w_i f'_j(u)  \frac{\partial^2 \phi_\delta}{\partial \xi_i \partial \xi_k}(\mathbf{w}) \frac{\partial w_k}{\partial x_j} 
\\[5pt]
&=\frac{\delta^{2}}{\left(|\mathbf{w}|^{2}+\delta^{2}\right)^{3 / 2}}\Big(\varepsilon|\nabla \mathbf{w}|^{2} - w_i f'_j(u) \frac{\partial w_i}{\partial x_j} \Big)
\\[5pt]
& \geq-\frac{1}{4 \varepsilon}\left|f^{\prime}(u)\right|^{2} \frac{\delta^{2}|\mathbf{w}|^{2}}{\left(|\mathbf{w}|^{2}+\delta^{2}\right)^{3 / 2}} \geq -\frac{\delta \, L_f^{2} }{4 \varepsilon}. 
\end{aligned}
$$
Therefore, sending $\delta \to 0^+$ it follows from the above estimates that 
\begin{equation}
\label{eq_estimativa_gradW}
\begin{aligned}
\int_{\Omega}|\mathbf{w}(t)| d x  &\leq \int_{\Omega}|\mathbf{w}(0)| d x + L_f  \int_{0}^{t} \!\! \int_{\Gamma}  |\nabla u_b| \, dr dt  + \int_{0}^{t} \!\! \int_{\Gamma}  |D^2 u_b| \, dr dt 
\\[5pt]
&+C_{n,s}\int_0^t \!\! \int_{\Omega} \! \int_\Gamma \frac{|A(u(\tau,x)) - A(u_b(\tau,r))|}{|r-x|^{n+2s}} \, dr dx d\tau, 
\end{aligned}
\end{equation}
where we have used Lemma \ref{lem:estimati_positivo}. Consequently, the proof of the theorem is complete, since from item $(i)$, 
$$
C_{n,s} \! \! \int_0^t \! \! \int_{\Omega} \! \int_{\Gamma}\frac{|A(u(\tau,x))-A(u_b(\tau,x))|}{|x-r|^{n+2s}}drdxd\tau
\leq  \int_\Omega \varpi(0) \ dx + \int_{0}^{t}\!\! \int_{\Omega} |g| d x d \tau. 
$$
Indeed, we observe that 
$$
\begin{aligned}
   |A(u(\cdot,x)) - A(u_b(\cdot,r))| &- |A(u_b(\cdot,x)) - A(u_b(\cdot,r))| 
\\[5pt]   
   &\leq |A(u(\cdot,x)) - A(u_b(\cdot,x))|, 
\end{aligned}
$$
 thus for each $t \in (0,T)$ 
$$
\begin{aligned}
\int_{\Omega} |\nabla u(t)| d x & \leq \int_{\Omega} |\nabla u_0| d x + \int_{0}^{t} \!\! \int_{\Gamma} \big(L_f |\nabla u_b| + |D^2 u_b| \big) dr dt + (a+b) |\Omega| 
\\[5pt]
&+ \int_{0}^{t}\!\! \int_{\Omega} \big(| \partial_t u_b| + L_f |\nabla u_b| + |\Delta u_b| + L_A |(-\Delta)^s_\Omega u_b| \big) \, dx dt
\\[5pt]
&+  L_A \! \! \int_0^t \! \! \int_{\Gamma} C_{n,s} \! \int_{\Omega}\frac{|u_b(\tau,x)-u_b(\tau,r)|}{|x-r|^{n+2s}} \, dx dr d\tau. 
\end{aligned}
$$
\end{proof}

Finally, we have the following proposition, from which we establish the $L^1-$type contraction \eqref{L1-contration}. 
Again, we drop the subscript $\ve$ in 
$u_\ve$, $u_{0,\ve}$, $u_{b,\ve}$, and $a_\ve$, $b_\ve$, see \eqref{abepsilon}. 

\begin{proposition}
\label{Th_contraction}
Let $u$ and $v$ be  solutions of the problem \eqref{eq:regular} 
with initial data $u_0$, $v_0$, and boundary condition $u_b$, $v_b$ respectively. 
Let $\xi_\ve$, $L$ be defined by \eqref{XI4TECHOS}. Then, for each $t \in (0,T)$,
\begin{equation}
\label{contration+}
\begin{aligned}
\int_{\Omega}|u(t)-v(t)|^{+}\xi_{\varepsilon} \, dx   
\leq&\,\int_{\Omega}|u_0-v_0|^{+}\xi_{\varepsilon} \, dx 
\\
&+(L_f+\varepsilon L)\int_0^t\!\!\int_{\Gamma}|u_b-v_b|^{+} \, drd\tau 
\\
&-\frac{C_{n,s}}{2}\int_{0}^t \left[|A(u)(\tau)-A(v)(\tau)|^+,\xi_\ve\right]_{W^{s,2}(\Omega)} d\tau.
\end{aligned}
\end{equation}
and
\begin{equation}
\label{contration+-}
\begin{aligned}
\int_{\Omega}|u(t)-v(t)|\xi_{\varepsilon} \, dx   
\leq&\,\int_{\Omega}|u_0-v_0|\xi_{\varepsilon} \, dx 
\\
&+(L_f+\varepsilon L)\int_0^t\!\!\int_{\Gamma}|u_b-v_b| \, drd\tau 
\\
&-\frac{C_{n,s}}{2}\int_{0}^t \left[|A(u)(\tau)-A(v)(\tau)|,\xi_\ve\right]_{W^{s,2}(\Omega)} d\tau.
\end{aligned}
\end{equation}
\end{proposition}
\begin{proof}
1. Let us denote $w= u-v$, and define $w_0:= u_0 - v_0$, $w_b:= u_b - v_b$. We are going to prove \eqref{contration+},
similarly we have \eqref{contration+-}. 
From equation $(\ref{eq:regular})_1$, we obtain
\begin{equation}
w_t+\dive(f(u)-f(v))+\varepsilon(-\Delta)w+(-\Delta)^s_{\Omega}(A(u)-A(v))=0.\label{eq:eq(u-v)}
\end{equation}
\medskip
Then, multiplying the equation \eqref{eq:eq(u-v)}
by $\varphi^\prime_\delta(w)\xi_\varepsilon$, where 
$$
\varphi_\delta(z):=\left\lbrace
\begin{aligned}
(z^2+\delta^2)^{1/2}-\delta,&\quad z\geq 0,
\\[5pt]
0,&\quad \, z\leq 0,
\end{aligned}
\right.
$$
and integrating over $(0,t)\times\Omega$ we have
\begin{equation}
\label{EQCONT}
\begin{aligned}
&\int_{\Omega} \varphi_{\delta}(w(t,x)) \, \xi_\varepsilon(x) \, d x 
-\int_{\Omega}\varphi_{\delta}(w_0(x)) \, \xi_\varepsilon(x) \, d x
\\[5pt]
&-\int_0^t \!\! \int_\Omega \varphi_{\delta}^{\prime \prime}(w(\tau,x)) \, \xi_\varepsilon(x) \, (f(u)-f(v)) \cdot \nabla w(\tau,x) \, d x d\tau 
\\[5pt]
&-\int_0^t \!\! \int_\Omega \varphi_{\delta}^{\prime}(w(\tau,x)) \ (f(u)-f(v))   \cdot \nabla \xi_\varepsilon(x) \, d x d\tau 
\\[5pt]
&+\int_{0}^t \!\! \int_{\Omega} \xi_\varepsilon(x)
\varphi_{\delta}^{\prime}(w) \, (-\Delta)^s_{\Omega}\big( A(u)-A(v)\big) \, dx d\tau 
\\[5pt]
&+\varepsilon \! \int_0^t \!\! \int_\Omega \nabla\varphi_{\delta}(w) \cdot \nabla \xi_\varepsilon \, dx d\tau 
+\varepsilon \! \int_0^t \!\! \int_\Omega |\nabla w(\tau,x)|^{2} \varphi_{\delta}^{\prime \prime}(w)\xi_\varepsilon \, d x d\tau= 0.
\end{aligned}
\end{equation}

2. Now, we proceed with the standard inequality, that is to say 
\[
\begin{aligned}
-\varphi_{\delta}^{\prime \prime}(w) \, (f(u)-f(v)) \cdot \nabla w&+\varepsilon|\nabla w|^{2} \varphi_{\delta}^{\prime \prime}(w)
 \\[5pt]
& \geq\left\{-L_f|u-v|\,|\nabla w|+\varepsilon|\nabla w|^{2}\right\} \varphi_{\delta}^{\prime \prime}(w) 
\\[5pt]
& \geq-\frac{1}{4 \varepsilon}(u-v)^2 L_f^2 \varphi_{\delta}^{\prime \prime}(w)
\geq-\frac{\delta \, L_f^2}{4 \varepsilon}, 
\end{aligned}
\]
where we have used that, $z^2\varphi_\delta^{\prime\prime}(z)< \delta$.
Analogously, since $|z|\,|\varphi_\delta^{\prime}(z)|\leq\varphi_{\delta}(z)$, we obtain
$$
\begin{aligned}
-\varphi_{\delta}^{\prime}(w(\tau,x)) \ (f(u)-f(v))   \cdot \nabla \xi_\varepsilon(x)
&\geq-|u-v|L_f|\varphi^{\prime}_{\delta}(w)|\,|\nabla\xi_\varepsilon|
\\[5pt]
&\geq-L_f\varphi_{\delta}(w)\,|\nabla\xi_\varepsilon|.
\end{aligned}
$$
Then, from \eqref{EQCONT} and the two inequalities above 
\begin{equation}
\label{EQCONT1}
\begin{aligned}
&\int_{\Omega} \varphi_\delta(w(t)) \xi_\varepsilon \ dx
+\int_{0}^t \!\! \int_{\Omega} \, \xi_\varepsilon
\varphi_{\delta}^{\prime}(w(\tau))  \,  (-\Delta)^s_{\Omega}\big( A(u)-A(v)\big) \, dx d\tau
\\[5pt]
&\quad  -\frac{\delta \, T \, L_f^2}{4\varepsilon}\int_{\Omega}\xi_\varepsilon \ dx
-L_f\int_0^t \!\! \int_{\Omega}\varphi_\delta(w) \, |\nabla \xi_\varepsilon|dxd\tau
\\[5pt]
&\quad +\varepsilon\int_0^t \!\! \int_{\Omega}\nabla\varphi_\delta(w) \cdot \nabla\xi_\varepsilon dxd\tau 
\leq \int_{\Omega}\varphi_{\delta}(w_0) \,\xi_{\varepsilon} \, dx. 
\end{aligned}
\end{equation}

3. Finally, from \eqref{EQCONT1} and taking $\beta= \varphi_{\delta}$ in \eqref{preliminar}, we have 
$$
\begin{aligned}
\int_{\Omega} \varphi_\delta(w(t)) \, \xi_\varepsilon \, dx&
+ \int_{0}^t \!\! \int_{\Omega}\xi_\varepsilon(x)
\varphi_{\delta}^{\prime}(w(\tau))  (-\Delta)^s_{\Omega}\big( A(u)-A(v)\big) \, dx d\tau
\\[5pt]
&\leq\int_{\Omega}\varphi_{\delta}(w_0)\xi_{\varepsilon} \, dx
+\frac{\delta T L_f^2}{4\varepsilon}\int_{\Omega} \xi_\varepsilon \, dx
\\[5pt]
&+(L_f+L\varepsilon) \int_0^t \!\! \int_{\Gamma} \varphi_\delta(w_b(\tau)) \, dr d\tau.
\end{aligned}
$$

 \underline {Claim:} For almost all $(t,x) \in Q_T$, 
 $$
     \sgn^+(u(t,x)-v(t,x)) (-\Delta)^s_{\Omega}\big( A(u)-A(v)\big) \geq (-\Delta)^s_{\Omega}\big| A(u)-A(v)\big|^+.  
 $$
Proof of Claim: The proof is similar to that one in Lemma \ref{lem:sgn}. 

\smallskip
Then, passing to the limit as $\delta \to 0^+$ and using  Theorem \ref{Thintpart}, we get \eqref{contration+}.
\end{proof}

\subsection{Proof of Main Theorem}
\label{proofmainthm}

Now we are ready to prove the main theorem of this paper.

\begin{proof}{(Main Theorem)} 1. Let $\{u_\ve\}$, $\ve \in (0,1)$, be the solutions of the initial boundary 
value problem \eqref{eq:regular}. Due to estimates \eqref{ESTIME} and \eqref{ESTX} 
we can show that $\{u_\ve\}$ (for a subsequence) is precompact in $C^0([0,T]; L^1(\Omega))$. Indeed, 
it follows from a routine argument an application of the Arzel\`a-Ascoli Theorem. 
Therefore, there exists a function $u \in C^0([0,T]; L^1(\Omega))$, such that, (to an appropriate subsequence) 
$u_\ve(t,x)$ converges to $u(t,x)$ as $\ve \to 0$ for a.e. $(t,x) \in Q_T$. Morevoer, applying \eqref{MAXPRINCIPLE2}
it follows that, $u \in L^\infty(Q_T)$. 

\medskip
2. Now, let us show that $u$ satisfies \eqref{eq:solfrac00} and \eqref{eq:solfrac},
that is to say, $u$ is a week solution to the initial-boundary value problem \eqref{FTPME}. 
First, to show \eqref{eq:solfrac00} we multiply equation $(\ref{eq:regular})_1$ by $A(u_\ve)$, integrate in $\Omega$
and after some partial integration, we obtain for each $t \in (0,T)$, 
$$
\begin{aligned}
  &\ve \int_\Omega A^\prime(u_\ve(t)) \, |\nabla u_\ve(t)|^2 \, dx 
  + \frac{C_{n,s}}{2}\iint_{\Omega\times\Omega} \dfrac{|A(u_\ve(t,x))-A(u_\ve(t,y))|^2}{|x-y|^{n+2s}}dx\,dy
  \\[5pt] 
  &\leq  \int_\Omega |A(u_\ve(t))| \, |\partial_t u_\ve(t)| \, dx  
  + \int_\Omega |f(u_\ve(t))| \, A^\prime(u_\ve(t)) \, |\nabla u_\ve(t)| \, dx
  \\[5pt]
  &\leq C \big( L_A \,  C_1 + L_f \, L_A \, C_2 \big), 
\end{aligned}  
$$
where $C$ is a positive constant, and we have used the estimates \eqref{ESTIME}, \eqref{ESTX}.
Therefore, neglecting the first (positive) term in the left hand side of the above inequality, 
we pass to the limit as $\ve \to 0^+$ to obtain 
$$
 \ess \!\!\! \sup_{t \in (0,T)}  [A(u(t)), A(u(t))]_{W^{s,2}(\Omega)} < \infty. 
$$
Consequently, since $u \in L^\infty(Q_T)$, that is $\iint_{Q_T} | A(u)|^2  dx dt< \infty$, from \eqref{G_inner} 
we have $A(u) \in L^2(0,T;H^s(\Omega))$ . Moreover, since $A(u_\ve)=A(u_{b,\ve})$ on $\Gamma_T$, it follows that 
$$
A(u_\ve)-A(u_{b,\ve})\in L^2(0,T;H^s_0(\Omega)).
$$
Further, $A(u_\ve)-A(u_{b,\ve})$ is uniformly bounded in $L^2(0,T;H^s_0(\Omega))$ w.r.t.
$\ve$. Thus we obtain
(along suitable subsequence) that, $A(u_\ve)-A(u_{b,\ve})$ converge weakly to $A(u)-A(\tilde{u}_{b})$ in $L^2(0,T;H^s_0(\Omega))$, 
where we have used the uniqueness of the limit. Hence $u$ satisfies \eqref{eq:solfrac00}. 

\medskip
In order to show \eqref{eq:solfrac}, let $\phi\in C^{\infty}_0((-\infty,T)\times\R^n)$ be a nonnegative test function, 
such that, $\sgn^{\pm}\left(A(u_b) - A(k)\right) \, \phi= 0$ a.e. on $\Gamma_T$ for each $k \in \R$. Moreover, for each
Kru\v zkov semi-entropy flux pair $(\eta_k^{\pm},q_k^{\pm})$, it follows from Theorem \ref{EXISTHM}, that 
\begin{equation}
\label{QSFinal}
\begin{aligned}
&-\iint_{Q_T} \Big\{ \ \eta_k^\pm(u_{\varepsilon})\phi_t+ {q}_k^\pm(u_{\ve}) \cdot \nabla \phi- \ve \, \eta_k^\pm(u_{\ve})(-\Delta)\phi\Big\} \xi_{\varepsilon} \, dxdt
\\[5pt]
&+ \frac{C_{n,s}}{2} \iint_{Q_T} \!\! \xi_\ve(x) \!\! \int_{\Omega}   \big(\eta_{A(k)}^{\pm}(A(u_\ve(x))-\eta_{A(k)}^{\pm}(A(u_\ve(y)))\big)  
\\[5pt]
&\hspace{100pt} \times \frac{(\phi(x)-\phi(y))}{|x-y|^{n+2s}}\, dydx dt \leq I_1^\ve + 2 \, I_2^\ve + \frac{C_{n,s}}{2} \,  I_3^\ve, 
\end{aligned}
\end{equation}
where
$$
\begin{aligned}
I_1^\ve&:= \int_{\Omega}\eta_k^\pm(u_{0,\varepsilon}) \, \phi(0) \xi_{\varepsilon} \, dx
+ (L_f + \varepsilon L) \int_{\Gamma_T} \eta_k^\pm(u_{b,\varepsilon}) \, \phi \, drdt,
\\[5pt]
I_2^\ve&:= \varepsilon \iint_{Q_T} \eta_k^\pm(u_\ve) \nabla\phi \cdot \nabla\xi_\varepsilon \, dxdt, 
\\[5pt]
I_3^\ve&:= - \! \! \iint_{Q_T} \!\! \phi(t,x) \!\! \int_{\Omega}   \big(\eta_{A(k)}^{\pm}(A(u_\ve(x))-\eta_{A(k)}^{\pm}(A(u_\ve(y))\big)\frac{(\xi_\ve(x)-\xi_\ve(y))}{|x-y|^{n+2s}}\, dydx dt. 
\end{aligned} 
$$
To follow, we pass to the limit as $\ve \to 0^+$ in \eqref{QSFinal}. Since $u_\ve(t,x)$ converges to $u(t,x)$ as $\ve \to 0$ for a.e. $(t,x) \in Q_T$,
from \eqref{preliminar1}, standard arguments and applying the Dominated Convergence Theorem, 
we can pass to the limit in the left hand side of \eqref{QSFinal}, also in the $I_1^\ve$ term and show that $I_2^\ve$ converges to zero. It remains to study the $I_3^\ve$
term. From \eqref{Apx-Lema11} it follows for each $\delta> 0$ fixed that  
$$
\begin{aligned}
&| I_3^\ve| \leq  | \int_0^T \!\! \int_{\Omega_\delta} \! \int_{\Omega_\delta} \phi(t,x)\frac{(\xi_\ve(x)-\xi_\ve(y))(\Phi(x)-\Phi(y))}{|x-y|^{n+2s}} \, dy dx dt|
\\[5pt]
&\qquad +| \int_0^T \!\! \int_{\Omega_\delta} \! \int_{\Omega_\delta^c}\Phi(y)\frac{\big(\phi(t,x)-\phi(t,y)\big)\big(\xi_\ve(x)-\xi_\ve(y)\big)}{|x-y|^{n+2s}} \, dx dy dt |
\\[5pt]
&\qquad + | \int_0^T \!\! \int_{\Omega_\delta} \! \int_{\Omega_\delta^c} \phi(t,y) \, \Phi(y)\frac{\big(\xi_\ve(x)-\xi_\ve(y)\big)}{|x-y|^{n+2s}} \, dx dy dt |=:  |I_{31}^\ve| + |I_{32}^\ve| + |I_{33}^\ve|
\end{aligned}
$$
with the obvious notation, where for each $t \in (0,T)$, $\Phi(\cdot)= \eta_{A(k)}^{\pm}(A(u_\ve(t,\cdot))$, 
(we keep this notation for simplicity), and we take $\Psi(\cdot)= \phi(t,\cdot)$. For $0< \ve< \delta$, 
we proceed a change of variables analogously in the proof of Lemma \ref{lem:Lap_lim_A_beta}. Then, 
denoting for convenience $\theta_x= r_x - \tau_x \ve \nu(r_x)$, analogously $\theta_y$,  
we have 
$$
\begin{aligned}
\lim_{\ve \to 0^+} |I_{31}^\ve|&= \lim_{\ve \to 0^+}| \int_0^T \!\! \ve^2 \!\! \int_0^\frac{\delta}{\ve} \!\! \int_0^\frac{\delta}{\ve} \!\! \int_{\Gamma} \! \int_{\Gamma} \phi(t,\theta_x) (\xi_\ve(\theta_x)-\xi_\ve(\theta_y))
\\[5pt]
& \hspace{70pt}  \times \frac{(\Phi(\theta_x)-\Phi(\theta_y))}{|\theta_x-\theta_y|^{n+2s}} \, dr_y dr_x d\tau_x d\tau_y dt + O(\delta)|= O(\delta), 
\\[5pt]
\lim_{\ve \to 0^+} |I_{33}^\ve|&= \lim_{\ve \to 0^+} | \! \int_0^T \!\! \ve \!\! \int_0^\frac{\delta}{\ve} \! \! \int_{\Omega_\delta^c} \!\! \phi(t,\theta_y) \, \Phi(\theta_y) \, \frac{\big(\xi_\ve(x)-\xi_\ve(\theta_y)\big)}{|x-\theta_y|^{n+2s}} \, dx dr_y d\tau_y dt + O(\delta)|
\\[5pt]
&= |\int_0^T \! \int_{\Omega_\delta^c} \int_\Gamma \phi(t,r) \, \, \frac{\eta_{A(k)}^{\pm}(A(u_b(t,r)))}{|x-r|^{n+2s}} \, dr dx dt + O(\delta)|= O(\delta), 
\end{aligned}
$$
where we have used, $\sgn^{\pm}\left(A(u_b(t,r)) - A(k)\right) \, \phi(t,r)= 0$ for a.e. $(t,r) \in \Gamma_T$.
Passing to the limit as $\delta \to 0^+$, the two above terms converge to zero and also  
$$
\begin{aligned}
\lim_{\delta \to 0^+} |I_{32}^\ve|&= \lim_{\delta \to 0^+} |\int_0^T \!\! \int_{\Omega_\delta} \! \int_{\Omega_\delta^c}\Phi(y)\frac{\big(\phi(t,x)-\phi(t,y)\big)\big(\xi_\ve(x)-\xi_\ve(y)\big)}{|x-y|^{n+2s}} \, dx dy dt|
\\[5pt]
&= \frac{ 2 \, \clg{N}_\sigma}{C_{n,s}}\,  |\!  \int_{\Gamma_T} \eta_{A(k)}^{\pm}(A(u_b(t,r))) \, \partial_\nu^\sigma \phi(t,r) \, dr dt|= 0, 
\end{aligned}
$$
where we have used (3.10) in the proof of Theorem 3.3 in \cite{Guan1}, p. 303, 
(see also p. 293 in that paper for the definitions of the sets $G_\delta$ and $G_\delta^\prime$).  
Consequently, we obtain equation \eqref{eq:solfrac} (in the limit) from the equation \eqref{QSFinal}. 

\medskip
3. Finally, we prove the $L^1-$type contraction \eqref{L1-contration}.
Under the conditions of Proposition \ref{Th_contraction}, 
we have from \eqref{contration+-} that 
$$
\begin{aligned}
\int_{\Omega}|u_\ve(t)-v_\ve(t)|\xi_{\varepsilon} \, dx   
&\leq\,\int_{\Omega}|u_{0,\ve}-v_{0,\ve}|\xi_{\varepsilon} \, dx 
\\[5pt]
& +(L_f+\varepsilon L)\int_0^t\!\!\int_{\Gamma}|u_{b,\ve}-v_{b,\ve}| \, dr d\tau  - J_1^\ve, 
\end{aligned}
$$
where 
$$
J_1^\ve = \frac{C_{n,s}}{2}\int_{0}^t \!\! \iint_{\Omega \times\Omega} \frac{\big(\Phi(\tau,x)-\Phi(\tau,y)\big)\big(\xi_\ve(x)-\xi_\ve(y)\big)}{|x-y|^{n+2s}} \, dy dx d\tau, 
$$
and we have denoted for convenience $\Phi(t,\cdot)= |A(u_\ve(t,\cdot)) - A(v_\ve(t,\cdot))|$.  
Moreover, applying \eqref{Apx-Lema12} we may write 
$$ 
\begin{aligned}
J_1^\ve &= \frac{C_{n,s}}{2}\int_{0}^t \!\!\int_{\Omega_\delta}\int_{\Omega_\delta}\frac{(\Phi(\tau,x)-\Phi(\tau,y))(\xi_\ve(x)-\xi_{\ve}(y))}{|x-y|^{n+2s}} \, dy dx d\tau
\\[5pt]
&+ C_{n,s}\int_{0}^t \!\! \int_{\Omega_\delta}\int_{\Omega_\delta^c}\frac{(\Phi(\tau,x)-\Phi(\tau,y))(\xi_\ve(x)-\xi_{\ve}(y))}{|x-y|^{n+2s}} \, dy dx d\tau. 
\end{aligned}
$$
Therefore, taking the limit as $\ve \to 0^+$ and proceeding as before in item 2, we obtain \eqref{L1-contration}. 
Recall that, we have here $A(\tilde{u}_b), A(\tilde{v}_b) \in L^2\left(0,T; H^{2s}(\Omega)\right)$, see also Remark \ref{TRACECOND}. 
Analogously, it follows from \eqref{contration+} that 
$$
\begin{aligned}
\int_{\Omega}|u(t,x)-v(t,x)|^+ \, dx
&\leq\,\int_{\Omega}|u_0(x)-v_0(x)|^+ \, dx
\\[5pt] 
&+ L_f \int_0^t\!\!\int_{\Gamma}|u_b(\tau,r)-v_b(\tau,r)|^+ drd\tau
\\[5pt]
&-\mathcal{N}_\sigma\int_0^t\int_\Gamma\partial_\nu^\sigma|A(u_b)-A(v_b)|^+(\tau,r) drd\tau.
\end{aligned}
$$
Consequently, since $|A(u_b)-A(v_b)|^+= \sgn^+(u_b-v_b) \big(A(u_b)-A(v_b)\big)$,
if $u_0 \leq v_0$ and $u_b \leq v_b$ almost everywhere, then 
$u(t,x) \leq v(t,x)$ for almost all $(t,x) \in Q_T$, which completes the proof. 

\end{proof}

\section{Appendix}
Here we fix some notation and background used in this paper.
Let $\Omega \subset \R^n$ be a bounded open set having smooth ($C^2$)
boundary $\Gamma$.
First, we define a level set function $h$, given for $\delta>0$ (sufficiently small) by 
$$
h(x):=
\left\{\begin{aligned}
\min (\operatorname{dist}(x, \Gamma), \delta) & \text { for } x \in \Omega, 
\\
-\min (\operatorname{dist}(x, \Gamma), \delta) & \text { for } x \in \mathbb{R}^{d} \backslash \Omega, 
\end{aligned}
\right.
$$
one can refer to \cite{Malek}. 
Therefore, the function $h(x)$ is Lipschitz continuous in $\mathbb{R}^{d}$, and 
smooth on the closure of $\left\{x \in \mathbb{R}^{d}:|h(x)|<\delta\right\}$.  
Also we have
\[
|\nabla h(x)|=\left\{\begin{array}{ll}
1 & \text { for } 0 \leq h(x)< \delta,
\\[5pt]
0 & \text { for } h(x)= \delta.
\end{array}\right.
\]

Now, for $\varepsilon>0$ we define conveniently $\xi_{\varepsilon}$ by
\begin{equation}
\label{XI4TECHOS}
\xi_{\varepsilon}(x):= 1-\exp \left(-\frac{L_f+\varepsilon L}{\varepsilon} h(x)\right), 
\end{equation}
where $L \equiv \sup _{0<h(x)<\delta}|\Delta h(x)|$ and $L_f>0$. The function $\xi_\ve(x) \equiv 0$ on $\Gamma$, 
and it satisfies the following weak differential inequality
\begin{equation}
\label{preliminar} 
   L_f \int_{\Omega}\left|\nabla \xi_{\varepsilon}\right| \beta 
   \leq \varepsilon \int_{\Omega} \nabla \xi_{\varepsilon} \cdot \nabla \beta+(L_f+L \varepsilon) \int_{\Gamma} \beta \, d r
\end{equation}
for each nonnegative test function $\beta \in C^\infty_c\left(\mathbb{R}^{d}\right)$, see \cite{Malek} p. 129. 
Also, we have 
\begin{equation}
\label{preliminar1}
\lim_{\ve\to0}\int_\Omega|\xi_\ve(x)-1|dx=0 \mbox{ and } 
\lim_{\ve\to0}\ve\int_\Omega|\nabla\xi_\ve(x)|dx=0.
\end{equation}
%
In the following, we have a more delicate result, with respect to $\xi_\ve(x)$,
which comes up from the application of the Regional Fractional Laplacian. 
\begin{lemma}
\label{Apx-Lema1} 
Let $\Psi$, $\Phi$ be non-negative smooth functions in $\R^n$, and $\xi_\ve$ as defined in \eqref{XI4TECHOS}. Then, we have
\begin{equation}
\label{Apx-Lema11} 
\begin{aligned}
& \int_{\Omega} \Psi(x) \int_{\Omega} \frac{(\xi_\ve(x)-\xi_\ve(y))(\Phi(x)-\Phi(y))}{|x-y|^{n+2s}} \, dy dx 
\\[5pt]
&\qquad \geq
\int_{\Omega_\delta}\int_{\Omega_\delta}\Psi(x)\frac{(\xi_\ve(x)-\xi_\ve(y))(\Phi(x)-\Phi(y))}{|x-y|^{n+2s}} \, dy dx
\\[5pt]
&\qquad -\int_{\Omega_\delta}\int_{\Omega_\delta^c}\Phi(y)\frac{\big(\Psi(x)-\Psi(y)\big)\big(\xi_\ve(x)-\xi_\ve(y)\big)}{|x-y|^{n+2s}} \, dx dy
\\[5pt]
&\qquad - \int_{\Omega_\delta}\int_{\Omega_\delta^c}\Psi(y)\Phi(y)\frac{\big(\xi_\ve(x)-\xi_\ve(y)\big)}{|x-y|^{n+2s}} \, dx dy, 
\end{aligned}
\end{equation}
where $\Omega_\delta$ is defined in \eqref{Omega_delta} and $\Omega_\delta^c= \Omega \setminus \Omega_\delta$.  
Similarly, we obtain
\begin{equation}
\label{Apx-Lema12} 
\begin{aligned}
\Big[\Phi\, ,\,\xi_\ve\Big]_{W^{s,2}(\Omega)}&= \iint_{\Omega\times\Omega}\frac{(\Phi(x)-\Phi(y))(\xi_\ve(x)-\xi_{\ve}(y))}{|x-y|^{n+2s}} \, dy dx
\\[5pt] 
&= \int_{\Omega_\delta}\int_{\Omega_\delta}\frac{(\Phi(x)-\Phi(y))(\xi_\ve(x)-\xi_{\ve}(y))}{|x-y|^{n+2s}} \, dy dx
\\[5pt]
&+ 2 \int_{\Omega_\delta}\int_{\Omega_\delta^c}\frac{(\Phi(x)-\Phi(y))(\xi_\ve(x)-\xi_{\ve}(y))}{|x-y|^{n+2s}} \, dy dx.
\end{aligned}
\end{equation}
\end{lemma}
\begin{proof}
We show \eqref{Apx-Lema11}, the proof of \eqref{Apx-Lema12} is similar. First, 
recall for each $\delta> 0$ (sufficiently small), the definition of
$\Omega_\delta$ and $\Omega_\delta^c$.  
Then, we have 
$$
\begin{aligned} 
& \int_{\Omega} \int_{\Omega}\Psi(x)\frac{(\xi_\ve(x)-\xi_\ve(y))(\Phi(x)-\Phi(y))}{|x-y|^{n+2s}} dx \, dy
\\[5pt]
&= \int_{\Omega_\delta}\int_{\Omega_\delta}\Psi(x)\frac{(\xi_\ve(x)-\xi_\ve(y))(\Phi(x)-\Phi(y))}{|x-y|^{n+2s}} dxdy + J_1 + J_2,
\end{aligned}
$$
where we have used that, $h(x)= h(y)= \delta$ for all $(x,y) \in \Omega_{\delta}^c \times \Omega_{\delta}^c$, 
which implies $\xi_\ve(x)= \xi_\ve(y)$ and
$$
\begin{aligned}
&J_1:= \int_{\Omega_\delta}\int_{\Omega_\delta^c}\Psi(x)\frac{(\xi_\ve(x)-\xi_\ve(y))(\Phi(x)-\Phi(y))}{|x-y|^{n+2s}} \, dx dy,
\\[5pt]
&
J_2:= \int_{\Omega_{\delta}^c}\int_{\Omega_{\delta}}\Psi(x)\frac{(\xi_\ve(x)-\xi_\ve(y))(\Phi(x)-\Phi(y))}{|x-y|^{n+2s}} \, dx dy. 
\end{aligned}
$$
Then, we have
$$
\begin{aligned}
J_1 \geq - \int_{\Omega_{\delta}}\int_{\Omega_{\delta}^c}\Psi(x)\Phi(y)\frac{(\xi_\ve(x)-\xi_\ve(y))}{|x-y|^{n+2s}}dxdy
\end{aligned}
$$
where we have used that $h(y)\leq \delta = h(x)$, thus $\xi_\ve(y)\leq\xi_\ve(x)$.
Similarly, due to $\xi_\ve(x)\leq\xi_\ve(y)$, it follows that
$$
\begin{aligned}
  J_2 &= - \int_{\Omega_{\delta}^c}\int_{\Omega_{\delta}}\Psi(x)\frac{(\xi_\ve(y)-\xi_\ve(x))(\Phi(x)-\Phi(y))}{|x-y|^{n+2s}}dxdy
\\[5pt]
&\geq - \int_{\Omega_{\delta}^c}\int_{\Omega_{\delta}}\Psi(x)\Phi(x)\frac{(\xi_\ve(y)-\xi_\ve(x))}{|x-y|^{n+2s}}dxdy
\\[5pt]
&=- \int_{\Omega_{\delta}}\int_{\Omega_{\delta}^c}\Psi(y)\Phi(y)\frac{(\xi_\ve(x)-\xi_\ve(y))}{|x-y|^{n+2s}}dxdy,
\end{aligned}
$$
from which the proof is complete.
\end{proof}

\subsection{Regional Fractional Parabolic Equation}

Now, let us consider an initial-boundary value problem, companion to \eqref{eq:regular}. More precisely, 
we study the existence of weak solutions of the following problem  
\begin{equation}
\begin{cases}
\partial_t u + \varepsilon(-\Delta)u+(-\Delta)_\Omega^sA(u+u_b)+\lambda u= f + f_b &\textit{in } Q_T,
\\[5pt]
u(0)=u_0 - u_b(0),&\textit{ in } \Omega,
\\[5pt]
u= 0, &\textit{ on } \Gamma_T,
\end{cases}\label{eq_reg_6.17}
\end{equation}
where $\lambda\geq 0$, $A^\prime(\cdot) \geq \varepsilon>0$, $u_0\in L^2(\Omega)$, $u_b \in H^1(Q_T)$, $f \in L^2(0,T; H^{-1}(\Omega))$, 
and 
$f_b= \Delta u_b - \partial_t u_b - \lambda u_b$.  
Then, we have the following 
\begin{theorem}
\label{Th_apen_Exis_2}
Under the above assumptions and denoting $\tilde{f}= f + f_b$, and 
$\tilde{u}_0= u_0 - u_b(0)$,
there exists a unique weak solution $u$ to \eqref{eq_reg_6.17}, that is to say, 
\begin{equation}
\left\lbrace
\begin{aligned}
&u\in L^2(0,T;H^1_0(\Omega)),\quad u_t\in L^2(0,T;H^{-1}(\Omega)), 
\\[5pt]
&\frac{d}{dt}(u(t),v) + B(u(t),v)=\left\langle \tilde{f}(t),v \right\rangle,\quad\textit{in }\mathcal{D}'(0,T), \quad (\forall v\in H^1_0(\Omega)), 
\\[5pt]
&u(0,x)= \tilde{u}_0, 
\end{aligned}
\right.\label{eq_th_ap_6.18}
\end{equation}
where $(\cdot,\cdot)$ denotes the scalar product in $L^2(\Omega)$, and $\left\langle\cdot,\cdot\right\rangle$ the paring between $H^{-1}(\Omega)$ and $H^1_0(\Omega)$. 
Moreover, the operator $B: H^1_0(\Omega) \times H^1_0(\Omega) \to \R$, is given by 
$$
\begin{aligned}
B(u,v):= \varepsilon & \int_{\Omega}\nabla u(x)\cdot\nabla v(x)\,dx
\\[5pt]
&+\frac{C_{n,s}}{2}
\iint_{\Omega\times\Omega}\!\!\frac{( A(u+u_b)(x)-A(u+u_b)(y))(v(x)-v(y))}{|x-y|^{n+2s}}dxdy
\\[5pt]
&+\lambda\int_{\Omega}u(x) v(x) \, dx. 
\end{aligned}
$$
\end{theorem}
\begin{proof}
1. First, we prove the existence result applying the Schauder fixed point argument. 
Let us consider $\Theta: \R^2 \to \R$, such that $\Theta(z_1,z_2)= \xi$, where  
$\xi= z$, when $z_1= z_2= z$, and otherwise
$$
   \xi= 
   \big(A^\prime\big)^{-1} \Big( \frac{A(z_1) - A(z_2)}{z_1 - z_2} \Big). 
$$
Then, given $w \in L^2(Q_T)$ we define $F(w) \equiv F_w$ by, 
$$
   F_w(t,x,y):= A^\prime \big(\Theta( (w+u_b)(t,x), (w+u_b)(t,y)) \big)
$$
for almost all $(t,x,y) \in (0,T) \times \Omega \times \Omega$, 
which satisfies, $\ve \leq F_w(t,x,y) \leq L_A$. Then, we consider the following auxiliary problem
\begin{equation}
\left\lbrace
\begin{aligned}
&\frac{d}{dt}(u(t),v)+B_w(u(t),v)= \left\langle g(t),v \right\rangle,\quad\textit{in }\mathcal{D}'(0,T), \quad (\forall v\in H^1_0(\Omega)), 
\\[5pt]
&u(0,x)= \tilde{u}_0
\end{aligned}
\right.\label{eq_th_ap_4_6.22}
\end{equation}
where 
$$
\begin{aligned}
B_w(u,v):= \varepsilon & \int_{\Omega}\nabla u(x)\cdot\nabla v(x)\,dx+\lambda\int_{\Omega}u(x)v(x)dx
\\[5pt]
&+\frac{C_{n,s}}{2}
\iint_{\Omega\times\Omega}\!\!
\frac{F_w(t,x,y)( u(x)-u(y))(v(x)-v(y))}{|x-y|^{n+2s}}dxdy,
\end{aligned}
$$
which is a bilinear mapping in $H^1_0(\Omega)$ and for a.a. $t \in (0,T)$
$$
\left\langle g(t),v \right\rangle= \left\langle \tilde{f}(t),v \right\rangle-\frac{C_{n,s}}{2}
\iint_{\Omega\times\Omega}\!\!
\frac{F_w(t,x,y)( u_b(x)-u_b(y))(v(x)-v(y))}{|x-y|^{n+2s}} \, dxdy.
$$
Moreover, we have that $B_w$ satisfies for all $u,v\in H^1_0(\Omega)$, 
$$
\begin{aligned}
|B_w(u,v)|
\leq 
\varepsilon & \|u\|_{H^1_0(\Omega)}\|v\|_{H^1_0(\Omega)}+\lambda\|u\|_{L^2(\Omega)}\|v\|_{L^2(\Omega)}
\\[5pt]
&+L_A\frac{C_{n,s}}{2}
\iint_{\Omega\times\Omega}\frac{|u(x)-u(y)||v(x)-v(y)|}{|x-y|^{n+2s}}dxdy
\\[5pt]
\leq C
 & \|u\|_{H^1_0(\Omega)}\|v\|_{H^1_0(\Omega)},
\end{aligned}
$$
where $C$ is a positive constant and also  
$$
\begin{aligned}
B_w(u,u)
= \varepsilon & \|u\|_{H^1_0(\Omega)}^2+\lambda\|u\|^2_{L^2(\Omega)}
\\[5pt]
&+\frac{C_{n,s}}{2}\int_{\Omega}\int_{\Omega}\frac{F_w(t,x,y)|u(x)-u(y)|^2}{|x-y|^{n+2s}} dxdy
\geq \varepsilon \|u\|_{H^1_0(\Omega)}^2. 
\end{aligned}
$$
Hence due to Theorem 11.7 in [Chipot], for each $w \in L^2(Q_T)$ fixed, 
the auxiliary problem \eqref{eq_th_ap_4_6.22} has a unique weak solution 
$$
u \in L^2(0,T;H^1_0(\Omega)),\quad \partial_t u \in L^2(0,T;H^{-1}(\Omega)), 
$$ 
and we have constructed the mapping
$$
   w \mapsto \mathcal{T}(w) \equiv u.
$$
Consequently, if $\mathcal{T}$ from $L^2(Q_T)$ into itself has a fixed point, then we are done with respect to the 
existence of solutions to \eqref{eq_th_ap_6.18}. To this end, it is enough to show that $\mathcal{T}(L^2(Q_T))$ is 
a relatively compact subset of the Banach space $L^2(Q_T)$, and also $\mathcal{T}$ is a continuous mapping 
with respect to the norm $\| \cdot \|_{L^2(Q_T)}$. 

\medskip
2. Now, for $R> 0$ to be chosen a posteriori, we denote by $B_R$ the following set
$$
B_R=\left\lbrace w\in L^2(0,T;L^2(\Omega)): \|w\|_{L^2(0,T;L^2(\Omega))}\leq R \right\rbrace.
$$
Then, let us show that $\mathcal{T}$ is a mapping from $B_R$ to itself. Indeed, 
taking $v= u(t)$ in \eqref{eq_th_ap_4_6.22}, we obtain
$$
\begin{aligned}
\frac{1}{2}\frac{d}{dt}\|u(t)\|^2_{L^2(\Omega)} + & \varepsilon  \, \|u(t)\|^2_{H^1_0(\Omega)} +
\lambda \|u(t)\|_{L^2(\Omega)}
\\[5pt]
&+\frac{C_{n,s}}{2}\int_{\Omega}\int_{\Omega}\frac{F(w)(t,x,y)|u(x)-u(y)|^2}{|x-y|^{n+2s}}dxdy=\left\langle g,u \right\rangle.
\end{aligned}
$$
Therefore, integrating on  $(0,T)$ and using the Young inequality, we have 
$$
\begin{aligned}
\|u(T)\|^2_{L^2(\Omega)}& + \varepsilon \int_0^T \|u(t)\|^2_{H^1_0(\Omega)} dt + 2 \lambda \int_0^T \|u(t)\|^2_{L^2(\Omega)} dt
\\[5pt]
& \leq \frac{1}{\ve}\|g\|_{L^2(0,T;H^{-1}(\Omega))} + \|\tilde{u}_0\|^2_{L^2(\Omega)}. 
\end{aligned}
$$
Thus there exists $R> 0$, which depends on $f$, $u_0$, $u_b$, $\lambda$, $\varepsilon$ and $C_{n,s}$, such that
\begin{equation}
\label{eq_th_ap_5_6.23}
\|u\|_{L^\infty(0,T;L^2(\Omega))}, 
\|u\|_{L^2(0,T;H^1_0(\Omega))}   \leq R,
\end{equation}
that is to say, $\mathcal{T}(w)\in B_R$ for each $w\in B_R$.
Moreover, from  \eqref{eq_th_ap_4_6.22}, we have the following estimate 
$$
\int_0^T\left\langle \partial_t u(t),v\right\rangle dt 
\leq C(R) \, \|v\|_{L^2(0,T;H^{1}_0(\Omega))},
$$
where $C(R)$ is a positive constant depending on $R$. 
Therefore, there exists $M>0$, such that
\begin{equation}
\|u_t\|_{L^2(0,T;H^{-1}(\Omega))}\leq M.\label{eq_th_ap_6_6.24}
\end{equation}
Thus, $u$ belongs to a bounded subset of $H^1(0,T;H^1_0(\Omega))$, which
is relatively compact in $L^2(Q_T)$. 
Therefore, in order to apply the Schauder fixed point argument in $B_R$, 
it remains to show that $\mathcal{T}$ is continuous on $B_R$. 

\medskip
Let $\{w_n\}$ be a sequence in $B_R$, such that, $w_n$ converges to $w$ in $B_R$ as $n\to \infty$,
and denote by $u_n= \mathcal{T}(w_n)$ the solution to \eqref{eq_th_ap_4_6.22} corresponding to $w_n$.
Due to \eqref{eq_th_ap_5_6.23} and \eqref{eq_th_ap_6_6.24} together with the Aubin-Lions' Theorem 
(along suitable subsequence), we have
\begin{equation}\label{eq_th_ap_7_6.25}
\begin{aligned}
w_n \to w &\quad \textit{a.e.} \quad (0,T) \times \Omega, 
\\[5pt]
u_n \to u &\quad \textit{in} \quad L^2(0,T;L^2(\Omega)),
\\[5pt]
\nabla u_n \rightharpoonup \nabla u &\quad \textit{in} \quad L^2(0,T;L^2(\Omega)), 
\\[5pt]
 \partial_t u_n \rightharpoonup  \partial_t u &\quad \textit{in} \quad L^2(0,T;H^{-1}(\Omega)). 
\end{aligned}
\end{equation}

\medskip
On the other hand, from \eqref{eq_th_ap_4_6.22}, we obtain
$$
\begin{aligned}
&\int_0^T \left\langle \partial_t u_n(t), v \right\rangle \, \zeta(t) \, dt + \varepsilon \int_0^T \!\!\int_{\Omega} \zeta(t) \, \nabla u_n(t) \cdot \nabla v \, dxdt
\\[5pt]
&+\lambda\int_0^T\!\!\int_{\Omega} \zeta(t) \, u_n(t) \, v \, dxdt
+\frac{C_{n,s}}{2}\int_0^T\!\!\!
\iint_{\Omega\times\Omega} \zeta(t) F(w_n) \, (u_n(t,x)-u_n(t,y))
\\[5pt]
& \qquad  \times \frac{(v(x)-v(y))}{|x-y|^{n+2s}} \, dxdydt
=\int_0^T\left\langle \tilde{f}(t),v \right\rangle \zeta(t) \, dt
\\[5pt]
&\qquad -\frac{C_{n,s}}{2}\int_0^T\!\!\!
\iint_{\Omega\times\Omega}\!\!\! \zeta(t) \, F(w_n) \, (u_b(t,x)-u_b(t,y)) \frac{(v(x)-v(y))}{|x-y|^{n+2s}} \, dxdydt
\end{aligned}
$$
for each $\zeta \in C^\infty_0(0,T)$. 
Then, from the convergences in \eqref{eq_th_ap_7_6.25}, the continuity regularity of $F$ and applying the Dominated Convergence Theorem,
we pass to the limit in the above equation as $n \to \infty$, hence $u= \mathcal{T}(w)$. 
Moreover, a routine argument shows that ,
$u(0)= \tilde{u}_0$. Therefore, the continuity of $\mathcal{T}$ is established.  
Consequently, by the Schauder fixed point argument, there exists $u \in B_R$ which satisfies \eqref{eq_th_ap_5_6.23} and \eqref{eq_th_ap_6_6.24}, such that $\mathcal{T}(u)=u$. 
Then, from \eqref{eq_th_ap_4_6.22} together with the definition of $F$, we obtain a weak solution to \eqref{eq_th_ap_6.18}. 

\medskip
3. Finally, we show the uniqueness, which follows from a standard routine. 
Indeed, we observe that equation \eqref{eq_reg_6.17} is equivalent to the problem
\begin{equation}
\begin{cases}
\partial_t u + \varepsilon(-\Delta)u+(-\Delta)_\Omega^sA(u)+\lambda u= f&\textit{in } Q_T,
\\[5pt]
u(0)=u_0 ,&\textit{ in } \Omega,
\\[5pt]
u= u_b, &\textit{ on } \Gamma_T,.
\end{cases}\label{eq_uniq_1}
\end{equation}
Then, let $u,v$ be two solutions to \eqref{eq_uniq_1} with the same data, 
and define $w=u-v$. Thus from $\eqref{eq_uniq_1}_1$, we obtain
\begin{equation}
\partial_t w + \varepsilon(-\Delta)w+(-\Delta)_\Omega^s\big(A(u)-A(v)\big)+\lambda w=0,\label{eq_uniq_2}
\end{equation}
in the distribution sense. Also we have $w(0)=0$ in $\Omega$ and $w=0$ on $\Gamma_T$. 
Then, we use $\eta(w)= |w|^+$ as a test function in \eqref{eq_uniq_2}, and from routine argument we obtain
for a.a. $t \in (0,T)$, 
$$
   \frac{d}{dt}\int_{\Omega}|w(t)|^+ \, dx 
   +\varepsilon \int_{\Omega}\eta''(w(t))|\nabla w(t)|^2 \, dx
   +\frac{C_{n,s}}{2}[A(u)-A(v), \eta^\prime(w)]_{W^{s,2}(\Omega)} \leq 0.
$$
Since the third term in the above equation is non-negatives, (see Lemma \ref{lem:sgn}), similarly the second one,  
it follows that 
$$
\int_{\Omega}|w(t,x)|^+dx \leq 0
$$
for a.a. $t \in (0,T)$. Therefore, we have $u\leq v$ almost everywhere. Exchanging the roles of $u$ and $v$, it leads to the uniqueness of solution.
\end{proof}

\section*{Data availability statement}
Data sharing is not applicable to this article as no data sets were generated or analysed during the current study.

 \section*{Conflict of Interest}
 Author Wladimir Neves has received research grants from CNPq
through the grant  308064/2019-4, and also by FAPERJ 
(Cientista do Nosso Estado) through the grant E-26/201.139/2021.


\end{document}